\documentclass[10pt]{amsart}
\usepackage[]{latexsym,amssymb,amsmath,amsfonts, amsthm}
\usepackage[all,cmtip,ps]{xy}
\usepackage{enumitem}

\def\N{{\mathbb N}}

\def\dualita#1#2{\mathrel{
                 \mathop{\vcenter{
                 \offinterlineskip
                 \hbox to 0.6truecm{\rightarrowfill}
                 \hbox to 0.6truecm{\leftarrowfill}}}%
                 \limits_{#2}^{#1}}}

\newtheorem{theorem}{Theorem}[section]

\newtheorem{corollary}[theorem]{Corollary}

\newtheorem{definition}[theorem]{Definition}
\newtheorem{example}[theorem]{Example}

\newtheorem{lemma}[theorem]{Lemma}

\newtheorem{proposition}[theorem]{Proposition}
\theoremstyle{remark}
\newtheorem{remark}[theorem]{Remark}

\newcommand{\grcong}{\cong_{gr}}

	\newcommand{\ZZ}{\mathbb{Z}}

	\newcommand{\NN}{\mathbb{N}}

\begin{document}

\title{Morita equivalence for graded rings}
\author{Gene Abrams}
\author{Efren Ruiz}
\author{Mark Tomforde}

\address{Department of Mathematics \\ University of Colorado \\ Colorado Springs, CO 80918-3733 \\USA}
\email{abrams@math.uccs.edu}

\address{Department of Mathematics\\University of Hawaii, Hilo\\200 W. Kawili St.\\ Hilo, Hawaii\\ 96720-4091 USA}
\email{ruize@hawaii.edu}

\address{Department of Mathematics \\ University of Colorado \\ Colorado Springs, CO 80918-3733 \\USA}
\email{mtomford@uccs.edu}

\date{\today}

\thanks{This work was supported by grants from the Simons Foundation (\#527708 to Mark Tomforde and \#567380 to Efren Ruiz)}


\begin{abstract}
The classical Morita Theorem for rings established the equivalence of three statements, involving categorical equivalences, isomorphisms between corners of finite matrix rings, and bimodule homomorphisms.  A fourth equivalent statement (established later) involves an isomorphism between infinite matrix rings.    In our main result, we establish the equivalence of analogous statements involving graded categorical equivalences, graded isomorphisms between corners of finite matrix rings, graded bimodule homomorphisms, and graded isomorphisms between infinite matrix rings.     
\end{abstract}

\maketitle

\section{Introduction}

Two (unital, associative) rings $R$ and $S$ are called {\it Morita equivalent}  when the categories of right modules $Mod-R$ and $Mod-S$ are equivalent.\footnote{We will work with categories of right modules here; but the categories of right modules over $R$ and $S$ are equivalent if and only if the categories of left modules are equivalent, so there is no loss of generality in doing so.}    In his foundational 1958 paper  \cite[Section 3]{Morita}, Morita proved that the following three conditions are equivalent:

\smallskip

\noindent {\bf The Original Morita Theorem}:

\begin{itemize}

\item[(M1)]  $R$ and $S$ are Morita equivalent (i.e., the categories $Mod-R$ and $Mod-S$ are equivalent).
 
\item[(M2)]  There exist $n\in \N$ and an idempotent $e \in {\rm M}_n(S)$ that is full in ${\rm M}_n(S)$ and for which the rings $R$ and $ e{\rm M}_n(S)e$ are isomorphic.
 
\item[(M3)] There exist an $R-S$-bimodule $P$ and an $S-R$-bimodule $Q$ and appropriate surjective bimodule homomorphisms  $P \otimes_S Q \to R$ and $Q\otimes_R P \to S$.   

\end{itemize}
 
\noindent For any ring $T$ we let ${\rm M}_\infty(T)$ denote the (nonunital) ring consisting  of those countably infinite square matrices over $T$  that contain  at most finitely many nonzero entries.  
 A fourth,  perhaps less well known condition, that is equivalent to those presented in The Original Morita Theorem,   was provided by Stephenson in his 1966 Ph.D. thesis \cite{Stephenson} (see also \cite{Abrams}).
 
\begin{itemize}
\item[(M4)]  The rings  ${\rm M}_\infty(R)$ and ${\rm M}_\infty(S)$ are isomorphic.        
\end{itemize}

\noindent We will refer to the equivalent statements (M1) through (M4) above as {\bf The Extended Morita 
 Theorem}.
 
In recent years there has been increased interest in graded rings, motivated in part by attempts to classify Leavitt path algebras and to use the natural $\mathbb{Z}$-grading on a Leavitt path algebra to study its structure.    As a result, it is natural to investigate Morita equivalence in the context of graded rings.  Suppose that $R$ is a ring that is graded by the group $\Gamma$, and let $Gr-R$ denote the category of graded right $R$-modules together with graded homomorphisms.   
   In \cite[Chapter 2]{RoozbehBook} the equivalence of the following two statements is established.

 \smallskip

\noindent {\bf The Graded Version of The Original Morita Theorem}:
 
\begin{itemize} 
\item[(GM1)]  The categories $Mod-R$ and $Mod-S$  are equivalent, via a functor that is compatible with a specific type of equivalence functor between the categories $Gr-R$ and $Gr-S$.

\item[(GM2)] There exist $n\in \N$ and an  idempotent $e \in {\rm M}_n(S)$ that is full in  ${\rm M}_n(S)$ and a sequence $(\gamma_m)_{1\leq m \leq n}$ in $ \Gamma$ for which the rings $R$ and $ e{\rm M}_n(S)[(\gamma_m)]e$ are graded isomorphic.
\end{itemize}
 
\noindent Although not explicitly included in \cite{RoozbehBook}, the following is easily seen to be equivalent to (GM1) and (GM2).   (See the comment following Theorem~\ref{GradedVersionExtMoritaThm}.)  

\begin{itemize} 
\item[(GM3)] There exist a graded $R-S$-bimodule $P$ and a graded $S-R$-bimodule $Q$ and appropriate surjective graded bimodule homomorphisms  $P \otimes_S Q \to R$ and $Q\otimes_R P \to S$.  
\end{itemize}

\noindent In Theorem~\ref{thm:gr-me-stablization} we establish the following analogue to (M4) of The Extended Morita Theorem in the graded setting and prove that it is equivalent to statements (GM1)--(GM3) of The Graded Version of the Original Morita Theorem.

\begin{itemize} 
\item[(GM4)] There exists a sequence $(\gamma_m)_{m \in \NN }$ in $\Gamma$ such that ${\rm M}_\infty(R)[(0)]$ is graded isomorphic to ${\rm M}_\infty(S) [ (\gamma_m) ]$.
\end{itemize}

\noindent 
An explanation of the notation used in (GM2) and (GM4) is given in Definition~\ref{sequencegradingdef} below. We will refer to the equivalent statements (GM1)--(GM4) above as {\bf The Graded Version of the Extended Morita Theorem}.  

Given a $\Gamma$-graded ring $R$, there are many ways to impose a $\Gamma$-grading on a matrix ring over $R$.  Arguably the most obvious is the {\it standard grading} in which the $\gamma$-component of the matrix ring is simply the matrix ring over the $\gamma$-component; i.e., $M_n(R)_\gamma = M_n(R_\gamma)$.   Furthermore, if $e$ is a homogeneous idempotent in $R$, then the corner ring $eRe$ directly inherits the $\Gamma$-grading from $R$, which we will refer to as the {\it standard} corner grading.    In the notation used in statements (GM2) and (GM4) above, the standard grading on the matrices arises when the sequences are constant (i.e., $\gamma_m = \gamma_{m'}$ for all $m,m'$).  

 It is natural to ask whether there is a collection of  statements about $\Gamma$-graded rings $R$ and $S$,  analogous to the  statements in The Graded Version of The Original Morita Theorem, in which the grading used on  each of the matrix rings is the standard grading.  We answer this question in the affirmative in our main result (Theorem \ref{MainTheorem}).     Specifically, we assume that the statement analogous to (GM4) holds, where the gradings on the matrix rings are standard; that is, we assume 
  $${\rm M}_\infty(R) \ \mbox{is graded isomorphic to} \  {\rm M}_\infty(S) \  \mbox{in the standard grading}.$$
  We then establish  the appropriate modifications  to statements (GM1), (GM2), and (GM3).     The power of Theorem \ref{MainTheorem} comes from the fact that there turns out to be a compatibility between the standard grading on matrix rings and  the existence of an appropriate category equivalence between the categories $Mod-R_0$ and $Mod-S_0$ of modules over the respective zero components of $R$ and $S$.      
 
 \smallskip

\noindent {\bf The homogeneously Graded Version of The Extended Morita Theorem} (See Theorem~\ref{MainTheorem}.)  
 
 \noindent Let $R$ and $S$ be $\Gamma$-graded rings.  Give each of the following rings of matrices over $R$ and $S$ the standard gradings.  Then the following are equivalent.

 \begin{itemize}
\item[(HG1)]   The categories $Mod-R$ and $Mod-S$  are equivalent as in (GM1), and the equivalence between $Gr-R$ and $Gr-S$ is,   in turn,  compatible with an equivalence of categories between $Mod-R_0$ and $Mod-S_0$.

 \item[(HG2)] There exist $n\in \N$ and a (homogeneous) idempotent $e \in {\rm M}_n(S)_0$ that is full in ${\rm M}_n(S)_0$   for which the rings $R$ and $ e{\rm M}_n(S)e$ are graded isomorphic.
 
\item[(HG3)] There exist a graded  $R-S$-bimodule $P$ and a graded  $S-R$-bimodule $Q$ and appropriate surjective graded bimodule homomorphisms  $P \otimes_S Q \to R$ and $Q\otimes_R P \to S$ that restrict to surjective bimodule homomorphisms of the germane $0$-components.

\item[(HG4)] ${\rm M}_\infty(R) $ is graded isomorphic to ${\rm M}_\infty(S)$.  
\end{itemize}

 \smallskip
 \noindent
 We call rings $R$ and $S$ satisfying condition (HG1) {\it homogeneously graded equivalent}.  
 \medskip
 
 We conclude the Introduction with a series of observations regarding our main results.  First, along the way to establishing  Theorem \ref{MainTheorem}, we give a new proof of the implication (M1) implies (M4) in the Extended Morita Theorem;  see Remark \ref{bottomupremark} below.

Second, we note that statement (GM2) does not yield statement (HG2) simply by imposing the standard grading on $ e{\rm M}_n(S)e$. As it turns out, we must also require that the idempotent in (HG2)   be full {\it in the zero component} ${\rm M}_n(S)_0$ of ${\rm M}_n(S)$.

Third, the notion of Morita equivalence is well known to functional analysts, and has been studied extensively for operator algebras.  In that context, Morita equivalence of the C$^*$-algebras $A$ and $B$ is defined by the existence of an imprimitivity Hilbert bimodule ${}_AX_B$.  The Brown-Green-Rieffel theorem \cite{BGR} states that $\sigma$-unital $C^*$-algebras $A$ and $B$ are Morita equivalent if and only if $A$ and $B$ are stably isomorphic (i.e., $A \otimes \mathcal{K} \cong B \otimes \mathcal{K}$, where $\mathcal{K}$ denotes the algebra of compact operators on a separable infinite-dimensional Hilbert space).  Since $\mathcal{K} = \overline{M_\infty(\mathbb{C})}$ and $A \otimes \mathcal{K} \cong \overline{M_\infty(A)}$, we see that $A \otimes \mathcal{K}$ is the analytic analogue of $M_\infty(A)$, and thus having $A$ stably isomorphic to $B$ (i.e., $A \otimes \mathcal{K} \cong B \otimes \mathcal{K}$) is the analytic analogue of having $M_\infty (A) \cong M_\infty (B)$.  So for  operator algebraists inquiring about corresponding ring-theoretic results, (HG3) and (HG4) are quite natural conditions.

Fourth, although one might expect to establish the graded versions of these results prior to establishing the homogeneously graded versions, we in fact proceed in the opposite order.  This is because the validity of one of the implications in the graded version will follow as a consequence  of the corresponding implication of the homogeneously graded version.  

 Fifth, we note that one of our motivating interests is when the rings $R$ and $S$ are Leavitt path algebras, and, more specifically, when the underlying graphs are finite with no sinks.  In this situation we prove that graded equivalence and homogeneously graded equivalence are identical (see Corollary \ref{strongLpa}.)  
 
 Sixth, a comment on terminology is in order.   Since we  are studying Morita Theory in the context of graded rings, it would not be inappropriate to call our investigation  ``Graded Morita Theory".   However,  that particular terminology has been used elsewhere in the literature (see, e.g., \cite{Boisen} and \cite{Haefner})   to denote something different from what we are studying in this article.  So we shall avoid using the (seemingly natural) phrase ``Graded Morita Theory".  
 
Seventh, we mention that applications of The Homogeneously Graded Version of the Extended Morita Theorem abound, and a number of these will be examined in a followup article.

Finally, the authors have found the book by R.~Hazrat \cite{RoozbehBook} to be an invaluable resource for many of the ideas we investigate in this paper, and readers will no doubt find this book useful and enlightening as well.  The authors also thank R.~Hazrat for personal communications and for many useful discussions while our article was being developed.

\section{The Algebraic Stabilization Theorem}

In this section we establish the equivalence of statements (HG2) and (HG4) of The Homogeneously Graded Version of The Extended Morita Theorem.   This result appears below  as Theorem~\ref{generalalgBGRthm}, and we refer to it as the {\it Algebraic Stabilization Theorem}.    

\begin{definition}\label{cornerdefinition}
Let $R$ be any associative (not necessarily unital) ring, and let $e = e^2$ be an idempotent in $R$.  The  ring $eRe := \{ exe \ | \ x\in R\}$ is called the {\rm corner ring of }$R$ {\rm generated by }$e$.   The corner ring $eRe$ is unital, with multiplicative identity $e$.   If $I$ is a two-sided ideal of $R$, then $eIe := \{ eje \ | \ j\in I\}$ is a two-sided ideal of $eRe$, and thus is a (not necessarily unital) ring in its own right. 

The idempotent $e$ of a ring $A$ is called {\rm full in} $A$ in case $AeA = A$;  that is, in case the two-sided ideal of $A$ generated by $e$ is all of $A$.  
\end{definition}

For two rings $A$ and $B$ graded by the abelian group $\Gamma$, we write $A \grcong B$ to mean there exists a $\Gamma$-graded ring isomorphism from $ A$ to $ B$.  For any ring $A$ we let ${\rm RCFM}(A)$ denote the ring of countably-infinite square matrices with entries in $A$ in which each row and each column contains at most finitely many nonzero entries.

\begin{definition}\label{standardgradingdef} Let $\Gamma$ be an abelian group, and let $A$ be a $\Gamma$-graded ring.  The \emph{standard matrix grading} on any matrix ring with entries in $A$  is  defined by setting $$({\rm M}_X(A))_\gamma = {\rm M}_X(A_\gamma).$$   (We abuse notation so that this definition  includes the ring ${\rm RCFM}(A)$ as well.)   Note that ${\rm M}_\infty(A)$ is a graded ideal of ${\rm RCFM}(A)$ in the standard matrix grading.
For an idempotent $e \in A_0$, the \emph{standard corner grading} on $eAe$ is given by setting $
(eAe)_\gamma = eA_\gamma e.$
\end{definition}

\begin{remark}\label{notfullinzerocomponent}   If $A$ is a unital graded ring and $e$ is an idempotent in $A_0$ for which $e$ is full in $A_0$  (i.e., if $A_0eA_0 = A_0$),   then $e$ is full in $A$ (since $1_A = 1_{A_0} \in A_0eA_0 \subseteq AeA$).   However, the converse need not be true:  an idempotent $e\in A_0$ that is full in $A$ need not be full in $A_0$.   As an easy example, consider $A = {\rm M}_2(K)$ for a field $K$, and grade $A$ by $\mathbb{Z}$ by setting $A_0$ to be the diagonal matrices, $A_1$ the upper corner matrices, and $A_{-1}$ the lower corner matrices.   Then $e = e_{1,1} \in A_0$ is full in $A$ (since $A$ is simple), but not full in $A_0$ (since $A_0 e A_0$ yields only the upper left corner matrices, and so does not contain $1_A$).    We note that $A$ is $\mathbb{Z}$-graded-isomorphic to the Leavitt path algebra $L_K(E)$, where $E$ is the graph  $ \bullet \longrightarrow \bullet$.  
\end{remark}

The main goal of this section is to prove the Algebraic Stabilization Theorem stated in Theorem \ref{generalalgBGRthm}.   The following result will provide nearly all of the heavy work towards establishing this goal.

\begin{theorem}\label{AlgebraicBGRspecialcase}
Let $\Gamma$ be an abelian group, let $A$ be a unital $\Gamma$-graded ring, and let $p \in A_0$ be an idempotent that is full in $A_0$; i.e.,  for which  $A_0 p A_0 = A_0$.  We grade  $pAp$ by the standard corner grading, and we grade ${\rm M}_\infty(A)$ and ${\rm M}_\infty(pAp)$ using the standard matrix grading.  Then ${\rm M}_\infty(A)$ is graded isomorphic to ${\rm M}_\infty(pAp)$. 
\end{theorem}

\bigskip

The proof of Theorem \ref{AlgebraicBGRspecialcase} will be given following Lemma \ref{lem:gr-corner-iso}; however, we outline the strategy here.  We construct idempotents $\{Q_n, P_n \ | \ n\in \mathbb{N}\}$ in the ring of row-finite and column-finite countably infinite matrices $\mathrm{RCFM}(A_0)$ together with graded homomorphisms 
   $$ \phi_n  \colon  \ Q_n {\rm M}_\infty(A) Q_n \to P_n {\rm M}_\infty(A) P_n \ \  \ \ \ \ \mbox{and} \ \ \ \  \ \  \psi_n \colon \ P_n {\rm M}_\infty(A) P_n \to Q_{n+1} {\rm M}_\infty(A) Q_{n+1}$$

\noindent
such that for all $n\in \mathbb{N}$ we have
\begin{enumerate}
    \item $Q_{n}Q_{n+1} = Q_n = Q_{n+1} Q_n$,
    
    \item $P_n P_{n+1} = P_n = P_{n+1} P_n$,
    
    \item $\bigcup_n P_n {\rm M}_\infty(A) P_n = {\rm M}_\infty( p A p )$,
    
    \item $\bigcup_n Q_n {\rm M}_\infty(A) Q_n = {\rm M}_\infty(A)$, and 
    
    \item the diagram
    
    \[
    \xymatrix{
    Q_n {\rm M}_\infty(A) Q_n \ar@{^{(}->}[rr]^{i_n} \ar[d]_-{\phi_n} & &  Q_{n+1} {\rm M}_\infty(A) Q_{n+1} \ar[d]^-{\phi_{n+1}}  \\
     P_{n} {\rm M}_\infty(A) P_n \ar[urr]^-{\psi_n} \ar@{^{(}->}[rr]^{j_n} & & P_{n+1} {\rm M}_\infty(A) P_{n+1}
    }
    \]
    commutes.   
\end{enumerate}
Conditions (1)--(4) imply
\[
{\rm M}_\infty(A) \grcong \varinjlim ( Q_n {\rm M}_\infty(A) Q_n , i_n ) \qquad  \text{and} \qquad {\rm M}_\infty(pAp) \grcong \varinjlim ( P_n {\rm M}_\infty(A) P_n , j_n )
\]
and the commutativity of the diagram in (5) implies that these direct limits are graded isomorphic.

\begin{definition} Let $\{ e_{i,j} \}$ be a system of matrix units for ${\rm M}_\infty( \ZZ )$.  For $a \in A$, we let $a \otimes e_{i, j}$ denote the element in ${\rm M}_n(A)$ or ${\rm M}_\infty(A)$ with $a$ in the $(i,j)$th entry and zero in all other entries.  For any infinite subset $S$ of $\NN$, $p \otimes 1_S := \sum_{ i \in S } p \otimes e_{i, i}$ is an element of $\mathrm{RCFM}(R)$ since $p\otimes 1_S$ is a countably infinite matrix with $p$ in the $(i,i)$ entry for all $i \in S$ and zero in all other entries.  

\noindent For idempotents $e$ and $f$ in a ring $A$, we write $e \leq f$ to mean $ef=fe=e$ (equivalently, that $e \in fAf$).  
\end{definition}

Note that if $u, v \in A$ such that $uvu=u$ and $vuv=v$, then $uv$ and $vu$ are easily seen to be idempotents in $A$.

\begin{lemma}\label{lem:sub-proj1}
Let $A$ be a ring.  Suppose $p$ and $e$ are idempotents in $A$ such that $uv = p$, $pu = u$, and $vu \in eAe$ for some $u, v \in A$.  Then there are $w,z \in A$ such that $wzw = w$, $zwz = z$, $wz = p$, and $zw \leq e$.  
\end{lemma}

\begin{proof}
Set $w = u$ and $z = vuv$.  Then $wz = uvuv=p^2=p$.  Since $pu=u$, we have $wzw = p u = u=w$, and $zwz = vuvuvuv=vuv=z$.  Therefore $zw$ is an idempotent.  Since $zw = vuvu \in eAe$, we conclude $zw \leq e$.  
\end{proof}
        
\begin{lemma}\label{lem:equiv-proj1}
Let $A$ be a unital ring, and let $p$ be an idempotent in $A$.  Also suppose $S$ and $T$ are subsets of $\NN$ with $|S| = |T| \in \N \cup \{\infty\}$.  Then there exist $u,v \in \mathrm{RCFM}(A)$ such that $uv=p \otimes 1_S$, $vu = p\otimes 1_T$, $uvu=u$, and $vuv=v$.  Moreover, if $|S|, |T| < \infty$, then $u, v$ can be chosen to be elements in ${\rm M}_\infty(A)$.
\end{lemma}

\begin{proof}
Let $m = |S| = |T|$.  Set $S = \{ s_i : i \in \NN, 1 \leq i \leq m  \}$ and $T= \{ t_i : i \in \NN, 1 \leq i \leq m \}$.  Define $u := \sum_{ i  = 1 }^m p \otimes e_{ s_i, t_i }$ and $v := \sum_{ i = 1}^m p \otimes e_{ t_i, s_i }$.  A straightforward computation shows that $u,v$ satisfy the desired conclusions.  The second statement is then clear.  
\end{proof}

\begin{lemma}\label{lem:sub-proj2}
Let $A$ be a unital ring, and let $p$ be a full idempotent in $A$.  Then there exists $m \in \NN$ and $w, z \in {\rm M}_\infty(A)$ such that $wzw=w$, $zwz=z$, $wz = 1_A \otimes e_{1,1}$, and $zw \leq \sum_{ i = 1}^m p \otimes e_{i,i}$.  Moreover, for any $j \in \NN$ and $S \subseteq \NN$ such that $|S| = m$, there exist $w',z' \in {\rm M}_\infty(A)$ such that $w'z'w'=w'$, $z'w'z'=z'$, $w'z' = 1_A \otimes e_{j, j}$, and $z'w' \leq p \otimes 1_{S}$.
\end{lemma}

\begin{proof}
Note that $1_A \in ApA$.  Therefore, $1_A = \sum_{ i = 1}^m x_i p y_i$ for some $x_i, y_i \in A$.  Set 
$$u := \left( \sum_{i=1}^m x_i \otimes e_{1, i} \right) \left( \sum_{i=1}^m p \otimes e_{i,i} \right) = \sum_{ i = 1}^m x_i p \otimes e_{1, i}$$
and  
$$v := \left( \sum_{i=1}^m p \otimes e_{i,i} \right) \left( \sum_{i=1}^m y_i \otimes e_{i, 1} \right) = \sum_{ i = 1}^m p y_i \otimes e_{i, 1}.$$
Then
\[
uv = \left( \sum_{ i = 1}^m x_i p \otimes e_{1, i}  \right) \left( \sum_{ i = 1}^m p y_i \otimes e_{i, 1} \right) =\left( \sum_{i=1}^m x_i p y_i \right) \otimes e_{1,1} = 1_A \otimes e_{1,1}.
\]
Note that $(1_A \otimes e_{1,1})u=u$ and $vu \in \left( \sum_{i = 1}^m, p \otimes e_{i.i} \right) {\rm M}_\infty(A) \left( \sum_{i = 1}^m p \otimes e_{i.i} \right)$.  By Lemma~\ref{lem:sub-proj1}, there are $w, z \in {\rm M}_\infty(A)$ with the desired properties.

Let $j \in \NN$ and let $S \subseteq \NN$ such that $|S| = m$.  By Lemma~\ref{lem:equiv-proj1}, there are $u_1, v_1$ and $u_2, v_2$ such that $u_iv_iu_i = u_i$, $v_iu_iv_i=v_i$, $u_1v_1 = 1_{A}\otimes e_{j,j}$, $v_1u_1 = 1_A \otimes e_{1,1}$, $u_2v_2 = \sum_{ i = 1}^m p \otimes e_{i,i}$ and $v_2 u_2 = p \otimes 1_{S}$.

Set $w' = u_1 w u_2$ and $z' = v_2 z v_1$.  Then

\begin{align*}
w' z' &= u_1 w u_2v_2 z v_1 = u_1 w 
\left( \sum_{ i =1}^m p \otimes e_{i,i} \right) z u_1  = 
   u_1 wzw
\left( \sum_{ i =1}^m p \otimes e_{i,i} \right) z u_1 \\
  & =  u_1 w zw z u_1 
     = u_1 (1_{A} \otimes e_{1,1}) v_1= u_1 v_1 u_1 v_1 = 1_{A} \otimes e_{j,j} ,      
\end{align*} 
so that 
$$w'z'w' = (1_{A} \otimes e_{j,j}) u_1 w  u_2 = u_1 w u_2 = w' ,  \ \ \ \  z' w'z' = v_2 z v_1 (1_{A} \otimes e_{j,j}) = v_2 z v_1 = z' ,  \ \ \mbox{and}$$
$$z'w' = v_2 (z v_1u_1 w )u_2 = v_2u_2 v_2 (z v_1u_1 w )u_2v_2u_2 \leq p \otimes 1_S.$$
This yields the second part of the lemma.  
\end{proof}

\begin{definition}
Two idempotents $p$ and $q$ in a ring $R$ are called  \emph{orthogonal} when $pq=qp=0$.
\end{definition}

\begin{lemma}\label{lem:sub-proj3}
Let $A$ be a unital ring, and let $p$ be a full  idempotent in $A$.  For any infinite subset $S$ of $\NN$, there are $w, z \in \mathrm{RCFM}(A)$ such that $wzw=w$, $zwz=z$, $wz = 1_A \otimes 1_S$ and $zw \leq p \otimes 1_S$.
\end{lemma}
\begin{proof}
Let $m \in \NN$ be the integer provided by Lemma~\ref{lem:sub-proj2}.  Let $\{ S_j \}_{j \in S }$ be disjoint subsets of $\NN$ such that $|S_j|= m$.  For each $j \in S$, by Lemma~\ref{lem:sub-proj2} there are $u_j, v_j \in {\rm M}_\infty(A)$ such that $u_jv_ju_j = u_j$, $v_j u_j v_j = v_j$, $u_j v_j = 1_A \otimes e_{j,j}$, and $v_j u_j \leq p \otimes 1_{S_j}$.  Set $u = \sum_{ j \in S} u_j$ and $v = \sum_{ j \in S } v_j$.  Since $u_j = u_j v_j u_j = (1_A \otimes e_{j,j}) u_j = u_j (p \otimes 1_{S_j})$ and $v_j = v_ju_j v_j = (p \otimes 1_{S_j}) v_j = v_j (1_A \otimes e_{j,j})$, $u$ and $v$ are elements in $\mathrm{RCFM}(A)$, and $\{ u_j v_j \}_{ j \in S }$ and $\{v_j u_j \}_{j \in S }$ are collections of mutually orthogonal idempotents.  Moreover:
\begin{align*}
uv &= \sum_{ j \in S} u_j v_j = \sum_{ j \in S} 1_A \otimes e_{j,j} = 1_A \otimes 1_S  \ ,\\
uvu &= (1_A \otimes 1_S) \sum_{ j \in S } (1_A \otimes e_{j,j}) u_j   = u \ , \\
vuv &= \left(\sum_{ j \in S } v_j (1_A \otimes e_{j,j}) \right) (1_A \otimes 1_S) = v \ , \ \ \ \mbox{and} \\
vu &= \sum_{ j \in S} v_j u_j  = \sum_{ j \in S} (p \otimes 1_{S_j}) v_j u_j (p \otimes 1_{S_j})  \\
    & \hspace{.5in} =( p \otimes 1_{X})\left(\sum_{ j \in S} (p \otimes 1_{S_j}) v_j u_j (p \otimes 1_{S_j})\right)( p \otimes 1_{X})   \leq p \otimes 1_{X}
\end{align*}
where $X = \bigsqcup_{ j \in S } S_j$.  By Lemma~\ref{lem:equiv-proj1}, there are $u',v' \in \mathrm{RCFM}(A)$ such that $u'v'u' = u'$, $v'u'v'=v'$, $u'v'= p \otimes 1_X$ and $v'u' = p \otimes 1_S$ since $|X|=|S|=\infty$.  

Set $w = uu'$ and $z=v'v$.  Then 
$$wz = u u' v' v = u (p \otimes 1_X) v = u vu (p \otimes 1_X) vuv  = u vuvu v  = uv = 1_A \otimes 1_S  \ ,  \ \ \ \mbox{so that} $$
$$wzw = (1_A \otimes 1_S) uu' = uvu u' = uu'=w \ , \ \    zwz = v'v (1_A \otimes 1_S) = v' vuv=v'v= z \ , \ \ \ \mbox{and} $$
$$ zw = v'v uu' = v'u'v' v u u'v'u' \leq p \otimes 1_S.$$ 
This proves the lemma.
\end{proof}

The following result plays a central role in our work.   The statement of the result, and a substantial portion of its proof, is an algebraic version of  \cite[Lemma~2.5]{Brown-PJM1977}.

\begin{theorem}[{\cite[Lemma~2.5 (Brown)]{Brown-PJM1977}}]\label{thm:brown}
Let $A$ be a unital ring, and let $p$ be a full idempotent in $A$.  Let $\NN = \bigsqcup_{i=1}^\infty N_i$ such that $|N_i|=\infty$ for each $i \in \N$.  Then there exist $u_i, v_i \in \mathrm{RCFM}(A)$ such that $u_i v_i u_i = u_i$, $v_i u_i v_i = v_i$, and the two sets $\{ u_j v_j \}_{ j \in \NN }$ and $\{v_j u_j \}_{j \in \NN }$ consist of  of mutually orthogonal idempotents.  In addition:  
  $$\sum_{ i = 1}^{2n-1} u_i v_i = \sum_{j=1}^n 1_A\otimes 1_{N_j}, \ \ \ 
\sum_{ i = 1}^{2n} u_i v_i \leq \sum_{j=1}^{n+1} 1_A\otimes 1_{N_j},$$
$$\sum_{ n = 1}^{2n-1} v_i u_i \leq  \sum_{j =1}^{n} p \otimes 1_{N_j}, \ \ \mbox{and}  \ \ 
\sum_{i=1}^{2n} v_i u_i = \sum_{j=1}^n p \otimes 1_{N_j}.$$
\end{theorem}        

\begin{proof}
We will inductively construct $u_i, v_i$ in $\mathrm{RCFM}(A)$.  By Lemma~\ref{lem:sub-proj3} there are $u_1, v_1 \in \mathrm{RCFM}(A)$ such that $u_1v_1u_1=u_1$, $v_1u_1v_1=v_1$, $u_1v_1 = 1_A \otimes 1_{N_1}$, and $v_1u_1 \leq p \otimes 1_{N_1}$.  By Lemma~\ref{lem:equiv-proj1} there are $w,z \in \mathrm{RCFM}(A)$ such that $wzw=w$, $zwz=z$, $wz = p \otimes 1_{N_2}$, and $zw = p \otimes 1_{N_1}$.  

Set $u_2 = w( p \otimes 1_{N_1} - v_1 u_1 )$ and $v_2 = (p \otimes 1_{N_1} - v_1u_1) z$.  Then 
\begin{align*}
u_2 v_2u_2 &= w( p \otimes 1_{N_1} - v_1 u_1 )zw( p \otimes 1_{N_1} - v_1 u_1 ) = w( p \otimes 1_{N_1} - v_1 u_1 ) = u_2 \ , \\
v_2 u_2 v_2 &= (p \otimes 1_{N_1} - v_1u_1) zw( p \otimes 1_{N_1} - v_1 u_1 )z = ( p \otimes 1_{N_1} - v_1 u_1 )z = v_2  \ , \\
u_2 v_2 &= w( p \otimes 1_{N_1} - v_1 u_1 )z = wzw( p \otimes 1_{N_1} - v_1 u_1 )zwz \leq 1_A \otimes 1_{N_2}  \ , \\
v_2u_2 &= (p \otimes 1_{N_1} - v_1u_1) z w( p \otimes 1_{N_1} - v_1 u_1 )  \\
    & \qquad \qquad = (p \otimes 1_{N_1} - v_1u_1)  p \otimes 1_{N_1} ( p \otimes 1_{N_1} - v_1 u_1 ) = p \otimes 1_{N_1} - v_1u_1 \ , \\
 u_1 v_1 u_2 v_2 &= (1_A \otimes 1_{N_1} )(1_A \otimes 1_{N_2} ) u_2 v_2 =0 \ , \\
u_2v_2 u_1v_1 &= u_2v_2 (1_A \otimes 1_{N_2} ) (1_A \otimes 1_{N_1} ) =0 \ , \\
v_1u_1 v_2u_2 &= v_1u_1 (p\otimes 1_{N_1} - v_1u_1) = 0 \ ,  \ \ \mbox{and}  \\
v_2u_2 v_1u_1 &= (p\otimes 1_{N_1} - v_1u_1)v_1u_1 =0.
\end{align*}
These equations imply that $u_1v_1$ and $u_2v_2$ are orthogonal idempotents, $v_1u_1$ and $v_2 u_2$ are orthogonal idempotents, $u_1 v_1 + u_2v_2 \leq 1_A \otimes 1_{N_1} + 1_A \otimes 1_{N_2}$, $v_1 u_1 \leq p \otimes 1_{N_1}$, and $v_1 u_1 + v_2u_2 = p \otimes 1_{N_1}$.  Suppose we have constructed $u_1, v_1, \ldots, u_{2n-2} , v_{2n-2}$ with the desired properties.  By Lemma~\ref{lem:sub-proj3} there are $u, v \in \mathrm{RCFM}(A)$ such that $uvu=u$, $vuv=v$, $uv= 1_A \otimes 1_{N_n}$ and $vu \leq p \otimes 1_{N_n}$.  Set
\begin{align*}
    u_{2n-1} &= \left( \sum_{j=1}^n ( 1_A \otimes 1_{N_j}) - \sum_{ k = 1}^{2n-2} u_kv_k \right) u \\
                &= \left( \sum_{j=1}^n ( 1_A \otimes 1_{N_j}) - \sum_{ k = 1}^{2n-2} u_kv_k \right) (1_A \otimes 1_{N_n} ) u \\
                &= ( 1_A \otimes 1_{N_n} - u_{2n-2} v_{2n-2} ) u
\end{align*}                
and   
\begin{align*}             
    v_{2n-1} &= v \left( \sum_{j=1}^n ( 1_A \otimes 1_{N_j}) - \sum_{ k = 1}^{2n-2} u_kv_k \right) \\
            &= v (1_A \otimes 1_{N_n} ) \left( \sum_{j=1}^n ( 1_A \otimes 1_{N_j}) - \sum_{ k = 1}^{2n-2} u_kv_k \right) \\
            &= v ( 1_A \otimes 1_{N_n} - u_{2n-2} v_{2n-2} ).
\end{align*}
Then
\begin{align*}
u_{2n-1} v_{2n-1} u_{2n-1}&= ( 1_A \otimes 1_{N_n} - u_{2n-2} v_{2n-2} ) u v ( 1_A \otimes 1_{N_n} - u_{2n-2} v_{2n-2} )u \\
&= ( 1_A \otimes 1_{N_n} - u_{2n-2} v_{2n-2} ) u = u_{2n-1} \ ,  \\
v_{2n-1} u_{2n-1} v_{2n-1} &=v ( 1_A \otimes 1_{N_n} - u_{2n-2} v_{2n-2} ) uv( 1_A \otimes 1_{N_n} - u_{2n-2} v_{2n-2}) \   \\
&= v ( 1_A \otimes 1_{N_n} - u_{2n-2} v_{2n-2} ) =v_{2n-1} \ ,  \ \ \ \mbox{and} \\
u_{2n-1} v_{2n-1} &= ( 1_A \otimes 1_{N_n} - u_{2n-2} v_{2n-2} ) u v ( 1_A \otimes 1_{N_n} - u_{2n-2} v_{2n-2} ) \\
&= 1_A \otimes 1_{N_n} - u_{2n-2} v_{2n-2}.
\end{align*}
Therefore, $u_{2n-1} v_{2n-1}$ is orthogonal to $u_k v_k$ for $1 \leq k \leq 2n-2$, and 
\[
\sum_{i=1}^{2n-1} u_i v_i = \sum_{i=1}^{2n-2} u_i v_i + 1_A \otimes 1_{N_n} - u_{2n-2}v_{2n-2} = \sum_{i=1}^n 1_A \otimes 1_{N_i}.
\]
Since
$
v_{2n-1}u_{2n-1} = vuv(1_A \otimes 1_{N_n} - u_{2n-2} v_{2n-2} ) uvu \leq p \otimes 1_{N_n}
$, 
and since $\sum_{k=1}^{2n-2} v_k u_k = \sum_{ i = 1}^{n-1} p \otimes 1_{N_j}$, we conclude $v_{2n-1} u_{2n-1}$ is orthogonal to $v_k u_k$ for $1 \leq k \leq 2n-2$, and 
\[
\sum_{j=1}^{2n-1} v_k u_k = \sum_{ j = 1}^{n-1} p \otimes 1_{N_j} + v_{2n-1} u_{2n-1} \leq \sum_{ j = 1}^n p \otimes 1_{N_j}.
\]

\noindent
By Lemma~\ref{lem:equiv-proj1} there are $w, z \in \mathrm{RCFM}(A)$ such that $wzw=w$, $zwz=z$, $wz = p \otimes 1_{N_{n+1}}$, and $zw = p \otimes 1_{N_n}$.  Set 
\begin{align*}
    u_{2n} &= w \left(\sum_{ i = 1}^n p \otimes 1_{N_i} - \sum_{k=1}^{2n-1} v_k u_k \right) = w (p \otimes 1_{N_n})\left(\sum_{ i = 1}^n p \otimes 1_{N_i} - \sum_{k=1}^{2n-1} v_k u_k \right)     \\
            &= w ( p \otimes 1_{N_n} - v_{2n-1} u_{2n-1})
 \end{align*}         
and            
\begin{align*}      
    v_{2n} &= \left(\sum_{ i = 1}^n p \otimes 1_{N_i} - \sum_{k=1}^{2n-1} v_k u_k \right) z = \left(\sum_{ i = 1}^n p \otimes 1_{N_i} - \sum_{k=1}^{2n-1} v_k u_k \right) (p \otimes 1_{N_n})z \\
            &= ( p \otimes 1_{N_n} - v_{2n-1} u_{2n-1})z.
\end{align*}
A similar computation to the one above shows that $u_{2n}v_{2n} u_{2n} = u_{2n}$, $v_{2n} u_{2n} v_{2n} = v_{2n}$,  $u_{2n}v_{2n} \leq p \otimes 1_{N_{n+1}}$, and $v_{2n} u_{2n} = p \otimes 1_{N_n} - v_{2n-1} u_{2n-1}$.  Since $\sum_{ i = 1}^{2n-1} u_i v_i = \sum_{ j = 1}^n 1_A \otimes 1_{N_j}$ is orthogonal to $p\otimes 1_{N_{n+1}}$, we conclude $u_{2n} v_{2n}$ is orthogonal to $u_{k}v_k$ for all $1 \leq k \leq 2n-1$ and $\sum_{ i =1}^{2n} u_i v_i \leq \sum_{ j = 1}^{n+1} 1_A \otimes 1_{N_j}$.  Another computation further shows that $v_{2n}u_{2n}$ is orthogonal to $v_k u_k $ for $1 \leq k \leq 2n-1$ and $\sum_{ i = 1}^{2n} v_i u_i = \sum_{ j = 1}^n p \otimes 1_{N_j}$, as desired.
\end{proof}

\begin{lemma}\label{lem:corner-iso}
Let $A$ be a ring, and let $p$ and $q$ be idempotents in $A$.  Suppose there exist $u, v \in A$ such that $uv = p$ and $vu \in qAq$.  Then the function $\phi \colon pAp \to qAq$ defined by $\phi(x)= v x u$ is a ring homomorphism.  Moreover, the following hold:  
\begin{enumerate}
\item  If $vu$ is an idempotent, then $\phi$ is an isomorphism from $pAp$ onto $vuAvu$. 

\item  If $I$ is a two-sided ideal of $A$, the restriction of $\phi$ to $pIp$ gives a homomorphism from $pIp$ into $qIq$ that is an isomorphism from $pIp$ onto $vu I vu$ when $vu$ is an idempotent.
\end{enumerate}
\end{lemma}

\begin{proof}
For $x \in pAp$, we have $\phi(x)=vxu = vpxpu = vuvxuvu = (vu)(vxu)(vu) \in qAq$.  It is clear that $\phi$ is an additive group homomorphism.  For $x, y \in pAp$, 
\[
\phi(x)\phi(y)=vxuvyu = vxpyu = vxyu = \phi(xy).
\]
Therefore, $\phi$ is a homomorphism.  

Suppose $e=vu$ is an idempotent.  Define $\psi \colon eAe \to pAp$ by $\psi(x) = uxv$.  A similar proof as above shows that $\psi$ is a homomorphism.  For $x \in pAp$ and $z \in eAe$, we have
\[
u(v x u) v = pxp=x \qquad \text{and} \qquad
v(uzv)u = e z e =z.
\]
Therefore $\phi$ and $\psi$ are inverse functions.

Note that if $x \in pIp$, then $\phi(x)=vxu \in qAq \cap I = qIq$.  Thus, $\phi \vert_{ p I p }$ is a homomorphism and if $vu$ is an idempotent, then $\psi \vert_{ vu I vu }$ is a homomorphism from $vuIvu$ to $pIp$ and is the inverse of $\phi \vert_{ p I p }$.  This establishes the last part of the lemma.
\end{proof}

\begin{lemma}\label{lem:gr-corner-iso}
Let $\Gamma$ be an abelian group and let $A$ be a $\Gamma$-graded ring.  If $I$ is a graded ideal of $A$ and $p \in A_0$ is an idempotent, then $e I_\gamma e = I \cap e A_\gamma e = I_\gamma \cap e I e$.  Consequently, $eIe$ is a $\Gamma$-graded ring with $(eIe)_\gamma = e I_\gamma e$.

Moreover, if $p, q \in A_0$ are idempotents, and $u, v \in A_0$ such that $uv= p$ and $vu \in q A q$, then the function $\phi \colon p I p \to q I q$ defined by $\phi(x)=vxu$ is a graded homomorphism.
\end{lemma}

\begin{proof}
Let $x \in e I_\gamma e$.  Since $I$ is a graded ideal, $I_\gamma = A_\gamma \cap I$.  Thus $x \in I \cap e A_\gamma e$.  Let $y \in I \cap e A_\gamma e$.  Then $y \in I \cap A_\gamma = I_\gamma$, and $y = eye \in I$.  Thus, $y \in I_\gamma \cap eIe$.  It is clear that $I_\gamma \cap e I e \subseteq e I_\gamma e$.  Hence $e I_\gamma e = I \cap e A_\gamma e = I_\gamma \cap e I e$.

By Lemma~\ref{lem:corner-iso}, $\phi$ is a ring homomorphism.  Suppose $x \in p I_\gamma p$.  Then $vxu \in I \cap A_\gamma = I_\gamma$. Since $vxu \in qIq$, it follows that $vxu \in  I_\gamma \cap q Iq  = q I_\gamma q$.  Thus $\phi$ is a graded homomorphism.  
\end{proof}

With the machinery of the above lemmas, we are now in position to give the proof of Theorem~\ref{AlgebraicBGRspecialcase}.  

\smallskip

\noindent \emph{Proof of Theorem~\ref{AlgebraicBGRspecialcase}}.
Let $\NN = \bigsqcup_{i=1}^\infty N_i$ such that $|N_i|=\infty$ for each $i\in \mathbb{N}$.  Since by hypothesis  $A_0 p A_0 = A_0$, we may apply Theorem~\ref{thm:brown} to the ring $A_0$ to obtain $u_i, v_i \in \mathrm{RCFM}(A_0)$ such that $u_i v_i u_i = u_i$, $v_i u_i v_i = v_i$, and $\{ u_j v_j \}_{ j \in S }$ and $\{v_j u_j \}_{j \in S }$ are collections of mutually orthogonal idempotents.  As well, 
$$\sum_{ i = 1}^{2n-1} u_i v_i = \sum_{j=1}^n 1_A\otimes 1_{N_j} , \ \  \sum_{ i = 1}^{2n} u_i v_i \leq \sum_{j=1}^{n+1} 1_A\otimes 1_{N_j} ,  \ \ \sum_{ n = 1}^{2n-1} v_i u_i \leq  \sum_{j =1}^{n} p \otimes 1_{N_j}, \ \   \mbox{and}  \ \ \sum_{i=1}^{2n} v_i u_i = \sum_{j=1}^n p \otimes 1_{N_j}.$$

\noindent
Set 
$$Q_n = \sum_{ j =1}^n 1_A \otimes 1_{N_j} \ \ \ \  \mbox{and} \ \ \ \  P_n = \sum_{j=1}^n p \otimes 1_{N_j}.$$  
Since $1_A \otimes 1_{N_j}$ and $p \otimes 1_{N_j}$ are idempotents in $\mathrm{RCFM}(A_0)$, it follows that $Q_n$ and $P_n$ are idempotents in $\mathrm{RCFM}(A_0)$.   Note that the inclusions of $Q_n {\rm M}_\infty(A) Q_n$ into $Q_{n+1} {\rm M}_\infty(A) Q_{n+1}$ and of $P_n {\rm M}_\infty(A) P_n$ into $P_{n+1} {\rm M}_\infty(A) P_{n+1}$ are graded homomorphisms and that
\[
\bigcup_{n =1}^\infty Q_n {\rm M}_\infty(A) Q_n = {\rm M}_\infty(A) \qquad \text{and} \qquad \bigcup_{n =1}^\infty P_n {\rm M}_\infty(A) P_n = {\rm M}_\infty(pAp).
\]
Let $i_n$ denote the embedding from $Q_n {\rm M}_\infty(A) Q_n$ into $Q_{n+1} {\rm M}_\infty(A) Q_{n+1}$, and let $j_n$ denote the embedding from $P_n {\rm M}_\infty(A) P_n$ into $P_{n+1} {\rm M}_\infty(A) P_{n+1}$.

Set $w_n = \sum_{ k = 1}^n u_k$ and $z_n = \sum_{k =1}^n v_k$. It is straightforward to check  that  the following hold:
\begin{align*}
    w_n z_n w_n &= w_n,  & z_n w_n z_n &= z_n,  \\
    w_{n+1} z_n  &= w_n z_n,  & w_n z_{n+1} &= w_n z_n,  \\
    z_{n+1} w_n &= z_n w_n,  & z_n w_{n+1} &= z_{n} w_n,  \\
    w_{2n-1}z_{2n-1} &= Q_n,  & z_{2n} w_{2n} &= P_n,  \\
    z_{2n-1} w_{2n-1} &\leq P_n,  \ \ \mbox{and} & w_{2n}z_{2n} &\leq Q_{n+1}.
\end{align*}

\noindent By Lemma~\ref{lem:gr-corner-iso}
$$\phi_n \colon Q_n {\rm M}_\infty(A) Q_n \to P_n {\rm M}_\infty(A) P_n \ \mbox{via} \ \  \phi_n (x) = z_{2n-1} x w_{2n-1}$$
and
$$\psi_n \colon P_n {\rm M}_\infty(A) P_n \to Q_{n+1} {\rm M}_\infty(A) Q_{n+1} \ \mbox{via} \ \ \psi_n(y) = w_{2n} y z_{2n}$$
 are graded homomorphisms.  We claim that $\psi_n \circ \phi_n = i_n$ and $\phi_{n+1} \circ \psi_n = j_n$.  Let $x \in Q_n {\rm M}_\infty(A)Q_n$ and $y \in P_n {\rm M}_\infty(A) P_n$. Then 
\begin{align*}
    \psi_n \circ \phi_n(x) &= w_{2n} z_{2n-1} x w_{2n-1} z_{2n} = w_{2n-1} z_{2n-1} x w_{2n-1} z_{2n-1} \\
    &= Q_n x Q_n = x = i_n(x)
\end{align*} 
while    
\begin{align*}    
    \phi_{n+1} \circ \psi_n (y) &= z_{2n+1} w_{2n} y z_{2n} w_{2n+1} = z_{2n} w_{2n} y z_{2n} w_{2n} \\ 
    &= P_n y P_n = y = j_n(y).
\end{align*}
Consequently, the families $\{ \phi_n \}_{n\in \N}$ and $\{ \psi_n \}_{n\in \N}$ induce a graded isomorphism between $\varinjlim ( Q_n {\rm M}_\infty(A) Q_n , i_n )$ and $\varinjlim ( P_n {\rm M}_\infty(A) P_n , j_n )$.  Therefore, ${\rm M}_\infty(A)$ and ${\rm M}_\infty(pAp)$ are graded isomorphic since $\varinjlim ( Q_n {\rm M}_\infty(A) Q_n , i_n )$ is graded isomorphic to ${\rm M}_\infty(A)$, and $\varinjlim ( P_n {\rm M}_\infty(A) P_n , j_n )$ is graded isomorphic to ${\rm M}_\infty(pAp)$.
\qed

\smallskip

The powerful but delicate result established in Theorem \ref{AlgebraicBGRspecialcase} can now be used to obtain the main result of this section.

\begin{theorem}\label{generalalgBGRthm}  (The Algebraic Stabilization Theorem)  \  Let  $A$ and $B$ be  unital rings graded by the abelian group $\Gamma$.   Impose the standard grading  on  any matrix ring over $A$ or $B$, and on any corner ring of a graded ring.  Then the following are equivalent:   

\begin{itemize}
\item[(HG2)]   There exists a positive integer $n$ and an idempotent $e\in {\rm M}_n(B)_0$  that  is full in ${\rm M}_n(B)_0$,   for which $A \grcong e{\rm M}_n(B)e$. 

\item[(HG4)]    ${\rm M}_\infty(A) \grcong {\rm M}_\infty(B)$. 
\end{itemize}
\end{theorem}

\begin{proof}   
Suppose (HG4) holds.  Let $\varphi:  {\rm M}_\infty(A) \to {\rm M}_\infty(B)$ be a graded isomorphism.   Because the grading is standard we have $A \grcong e_{1,1}{\rm M}_\infty(A) e_{1,1}$.     Let $e$ denote $\varphi(e_{1,1})$.   Then via the restriction of $\varphi$ to $e_{1,1}{\rm M}_\infty(A) e_{1,1}$,  $A \grcong e {\rm M}_\infty(B) e$.  But $e\in {\rm M}_n(B)$ for some $n\in \mathbb{N}$, so we get $e {\rm M}_\infty(B) e = e {\rm M}_n(B)e$, so that $$A \grcong e {\rm M}_n(B) e.$$   Because $e_{1,1}$ is full in $({\rm M}_\infty(A))_0 = {\rm M}_\infty(A_0)$, $e$ is necessarily full in  $({\rm M}_\infty(B))_0 = {\rm M}_\infty(B_0)$.   

We claim that $e$ is actually also full in ${\rm M}_n(B)_0$, which will complete this direction of the proof.    Note that if $S, T\in {\rm M}_\infty(B)_0$ then $e_nSe$ and $eTe_n$ are in ${\rm M}_n(B)_0$ where $e_n $ denotes $e_{1,1} + \cdots + e_{nn}$.   So let $X \in {\rm M}_n(B)_0$.  Then $e_nXe_n = X$.  So an equation $X = \sum_{i=1}^n S_i e T_i$ with $S, T\in {\rm M}_\infty(B)_0 = {\rm M}_\infty(B_0)$ (which exists by the fullness of $e$ in $({\rm M}_\infty(B))_0$) gives $X = e_nXe_n =  \sum_{i=1}^n (e_n S_i e) e (e T_i e_n)$, which gives the result.  Thus (HG2) holds.

\medskip

Conversely, suppose (HG2) holds.   Because we are utilizing the standard grading on all matrix rings, a  graded isomorphism   $A \grcong e {\rm M}_n(B) e$ induces a graded isomorphism (defined entrywise)  ${\rm M}_\infty(A) \grcong {\rm M}_\infty(e {\rm M}_n(B) e)$.  Since $e$ is full in ${\rm M}_n(B)_0$, Theorem  \ref{AlgebraicBGRspecialcase} gives that ${\rm M}_\infty(e{\rm M}_n(B)e) \grcong {\rm M}_\infty({\rm M}_n(B))$.   But because the grading on each of ${\rm M}_\infty$ and ${\rm M}_n$ is  the standard grading,  we conclude  ${\rm M}_\infty({\rm M}_n(B)) \grcong {\rm M}_\infty(B)$.   Thus
$${\rm M}_\infty(A) \grcong {\rm M}_\infty(e {\rm M}_n(B) e) \grcong {\rm M}_\infty({\rm M}_n(B)) \grcong {\rm M}_\infty(B),  $$    
\noindent
as desired.
\end{proof}

While Theorem \ref{generalalgBGRthm} is of interest in its own right, and provides the equivalence of two of the four statements that constitute Theorem \ref{MainTheorem},   there is an additional important consequence of Theorem \ref{generalalgBGRthm} that bears mentioning. If we choose $\Gamma$ to be the trivial group $\{0\}$, then ``$\Gamma$-graded" and ``ungraded" become synonymous.   In particular, in this situation the Algebraic Stabilization Theorem reduces to the equivalence of statements (M2) and (M4) in The Original Morita Theorem presented in the Introduction.

\begin{corollary}\label{ungradedcase}
Let $A$ and $B$ be unital rings.  Then ${\rm M}_\infty(A) \cong {\rm M}_\infty(B)$ if and only if $A$ and $B$ are Morita equivalent.  
\end{corollary}

\begin{remark}\label{bottomupremark}
It is compelling to note that the proof of Corollary \ref{ungradedcase} provided by Theorem \ref{generalalgBGRthm} follows a significantly different approach to the proofs of Corollary \ref{ungradedcase} that have appeared previously in the literature.  For any ring $R$ let $RFM(R)$ denote the ring of countably-square matrices over $R$ consisting of those matrices in which each row has at most finitely many nonzero elements.    Previously presented proofs have proceeded by showing that if $A$ and $B$ are Morita equivalent then there is an isomorphism between the rings $RFM(A)$ and $RFM(B)$; and then a tedious check confirms that this isomorphism restricts to an isomorphism between the nonunital subrings ${\rm M}_\infty(A)$ and ${\rm M}_\infty(B)$.   One may think of this as a ``top-down" approach to the isomorphism ${\rm M}_\infty(A) \cong {\rm M}_\infty(B)$.   However, the approach we provide here yields a ``bottom-up" approach,  achieved by interweaving various injective morphisms between finite-sized matrix rings over $A$ and $B$, respectively, and then taking a direct limit.     For additional information and historical perspective on these results, see  \cite{AS}.  
\end{remark}

We finish the section by presenting an additional  consequence of Theorem \ref{AlgebraicBGRspecialcase}, in the context of Leavitt path algebras.  (We refer the reader to \cite{AAS} for more information about this class of algebras.)     Suppose $E$ is a row-finite graph with no sinks (i.e., every vertex emits a finite and nonzero number of edges).  It is not hard to show (by invoking the (CK2) relation in $n$ stages)  that for any positive integer $n$ and for any vertex $w$ of $E$,  the element $w$ of the Leavitt path algebra $L_K(E)$ may be written as $w = \sum_{\mu \in S}\mu \mu^*$, where $S$ is the set of paths of length $n$ in $E$ having source vertex $w$.

\begin{definition}\label{primitivedef}
 Let $n$ be a positive integer.  A graph $E$ is called {\it primitive of length} $n$  if for any vertices $v,w \in E^0$, there exists a path from $v$ to $w$ having length exactly $n$.   We say that a graph $E$ is {\it primitive} if there is a positive integer $n$ such that $E$ is primitive of length $n$.  
 \end{definition}

\begin{corollary}\label{primitivecor}
Let $E$ be a primitive finite graph, and let $v$ be any vertex in $E$.   Impose the usual $\mathbb{Z}$-grading on $L_K(E)$, the induced grading on the corner $ vL_K(E)v$, and the standard grading on the matrix rings.   Then
$${\rm M}_\infty(L_K(E)) \grcong {\rm M}_\infty(vL_K(E)v).$$
\end{corollary}

\begin{proof}
By Theorem \ref{AlgebraicBGRspecialcase} it suffices to show that the idempotent $v$ of $L_K(E)_0$ is full in $L_K(E)_0$.  To do this we need only show that each vertex $w$ of $L_K(E)$ is in $L_K(E)_0 v L_K(E)_0$, as this will yield $ \sum_{w\in E^0} w = 1_{L_K(E)} \in L_K(E)_0 v L_K(E)_0$.     

Let $n$ be a positive integer such that $E$ is primitive of length $n$ (such a positive integer exists since $E$ is a primitive graph).  Let $\mu$ be any path in $E$ having ${\rm length}(\mu) = n$ and $s(\mu) = w$.   By hypothesis there exists a path $\nu$ in $E$  having ${\rm length}(\nu) = n$ for which $s(\nu) = v$ and $r(\nu) = r(\mu)$.     Let $x$ denote $\mu \nu^* \in L_K(E)$.   Since  ${\rm length}(\mu) = {\rm length}(\nu)$, $x \in L_K(E)_0$.  Since $s(\nu)=v$, $x^* = \nu\mu^* \in vL_K(E)_0$.  Since $x=xx^*x$ and since  $\mu \mu^* = xx^*$, we get that $x \in L_K(E)_0 v L_K(E)_0$ and hence, $\mu \mu^* \in  L_K(E)_0 v L_K(E)_0$.

We have shown that any path $\mu$ of $E$ with ${\rm length}(\mu) = n$ and $s(\mu) = w$, we have $\mu \mu^* \in L_K(E)_0 v L_K(E)_0$.   Since, by the remark in the paragraph before Definition~\ref{primitivedef}, $w$ can be written as  the sum of expressions of this form in $L_K(E)$, we have that $w\in L_K(E)_0vL_K(E)_0$, which yields the result. 
\end{proof}

\section{Homogeneously graded equivalence}  

In this section we use the Algebraic Stabilization Theorem  to establish a connection between the existence of a graded isomorphism ${\rm M}_\infty(A) \grcong {\rm M}_\infty(B)$ and the existence of a special type of equivalence functor between $Mod-A$ and $Mod-B$.    This will provide another piece of our main result, Theorem \ref{MainTheorem}.

For the reader's convenience, we provide a very brief primer of the key ideas we will need in this discussion.  A good resource for these concepts is  the book by Hazrat \cite{RoozbehBook}.   Throughout we let $\Gamma$ denote an abelian group.    If $A$ is a $\Gamma$-graded ring, then $Gr-A$ denotes the category whose objects are $\Gamma$-graded right $A$-modules and whose morphisms are $\Gamma$-graded homomorphisms.    For a graded right $A$-module $M$ and $\alpha \in \Gamma$, the \emph{$\alpha$-suspension of $M$}, denoted $M(\alpha)$, is the graded right $A$-module having $M(\alpha) = M$, with grading given by $M(\alpha)_\gamma = M_{\alpha + \gamma}$.    For $\alpha \in \Gamma$,    we let $\mathcal{T}_\alpha$ denote  the \emph{$\alpha$-suspension functor}  $\mathcal{T}_\alpha:  Gr-A \to Gr-A$  given by $M \mapsto M(\alpha)$ on objects and the identity on morphisms.

A functor $\phi: Gr-A \to Gr-B$ is called {\it graded} when $\phi \circ \mathcal{T}_\alpha = \mathcal{T}_\alpha \circ \phi $ for each $\alpha \in \Gamma$.     A graded functor $\phi : Gr-A \to Gr-B$ is a {\it graded equivalence} if there is a graded functor $\psi : Gr-B \to Gr-A$ such that $\phi$ and $\psi$ compose appropriately to the identity functors on each category.     If there is a graded equivalence between $Gr-A$ and $Gr-B$, we say $A$ and $B$ are {\it graded equivalent} or, more formally,  {\it graded Morita equivalent}.  

For any graded ring $A$, we let $U_A$ (or simply by $U$) denote the {\it forgetful functor} $U_A: Gr-A \to Mod-A$.  A functor $\phi^{\prime} : Mod-A \to Mod-B$ is called a {\it graded functor} if there is a graded functor $\phi: Gr-A \to Gr-B$ such that  $  U_B \circ \phi   =   \phi^{\prime} \circ U_A$ as functors from $Gr-A$ to $Mod-B$.  In this situation the functor $\phi$ is called an {\it associated graded functor} of $\phi^{\prime}$.    

A functor $\phi^\prime : Mod-A \to Mod-B$ is called a {\it graded equivalence} if it is both graded and an equivalence.

If $M_S$ is any right $S$-module and ${}_SN_T$ is a left $S$- right $T$- bimodule for which $N_T$ is $\Gamma$-graded, then we define a grading on  the right $T$-module  $M_S \otimes_S {}_SN_T$ by setting $(M_S \otimes_S {}_SN_T)_\gamma = M_S \otimes_S ( {}_SN_T)_\gamma$ for each $\gamma \in \Gamma$.

\begin{definition}\label{hgrdef}   {\rm We  call the $\Gamma$-graded rings $A$ and $B$ {\it homogeneously graded equivalent} in case there exists  a graded equivalence 
$$\psi: Gr-A \to Gr-B$$
for which there is an equivalence  
$$\eta: Mod-A_0 \to Mod-B_0$$
such that the diagram
\[   \xymatrix{Mod-A_0 \ar[d]_{- \otimes_{A_0}A}\ar[r]^\eta&Mod-B_0 \ar[d]^{- \otimes_{B_0}B}\\
Gr-A\ar[r]^\psi &Gr-B  }\]

\noindent
commutes  on objects of $Mod-A_0$ (up to isomorphism).  That is, $A$ and $B$ are called homogeneously graded equivalent in case  there is a category equivalence $\eta$ for which, for each object $M$ of $Mod-A_0$, there is an isomorphism  $\psi(M\otimes_{A_0}A) \grcong (\eta(M))\otimes_{B_0}B$ as objects of $Gr-B$.  In this situation we  write $A \approx_{hgr} B$. } 

\end{definition}

It is straightforward to show that ``homogeneously graded equivalent" is an equivalence relation on the collection of $\Gamma$-graded rings.

As described in the Introduction, we seek an analogue of The Extended Morita Theorem for graded rings.  With that in mind, we offer the following equivalent description of homogeneously graded equivalence, which follows directly from an application of  \cite[Theorem 2.3.7]{RoozbehBook}.

\begin{proposition}\label{equivalentversionhge}
Let $A$ and $B$ be $\Gamma$-graded rings.   Then $A$ and $B$ are homogeneously graded equivalent if and only if there exists a graded equivalence of categories $\psi': Mod-A \to Mod-B$ for which the associated functor $\psi: Gr-A \to Gr-B$ satisfies the property described in Definition \ref{hgrdef}.
\end{proposition}

Our goal is to establish a graded-module-theoretic framework that sheds light on the existence of a graded isomorphism between ${\rm M}_\infty(A)$ and $ {\rm M}_\infty(B)$.    We achieve this goal in Theorem~\ref{hgeTheorem}, and observe that the notion of homogeneously graded equivalence is precisely the required framework.

\begin{theorem}\label{hgeTheorem} Let $A$ and $B$ be $\Gamma$-graded rings.  Then  the following are equivalent.

\begin{itemize}
\item[(HG1)]      $A\approx_{hgr} B$.
\item[(HG4)]     ${\rm M}_\infty(A) \grcong {\rm M}_\infty(B)$ (in the standard grading).  
\end{itemize}

\end{theorem}

This section is devoted to proving Theorem~\ref{hgeTheorem} and examining some of its consequences.  For the reader's benefit, here is an outline our approach to proving Theorem~\ref{hgeTheorem}.  Similar to the way we proceeded   in the proof of the implication (HG1) $\Rightarrow$ (HG4) in  Theorem~\ref{generalalgBGRthm}, we start by analyzing the specific situation where $B = eAe$ for some idempotent $e$ of $A$.   We begin by showing that when $e$ is full in $A_0$, there is an appropriate diagram corresponding to $B = eAe$ and $B_0 = eA_0e$. Subsequently  we will show that the existence of an isomorphism ${\rm M}_\infty(A) \grcong {\rm M}_\infty(B)$ allows us to apply this specific result to deduce $A \approx_{hgr} B$.      

Conversely, we show $A\approx_{hgr} B$ implies that the existence of $n\in \mathbb{N}$ and an idempotent $e\in {\rm M}_n(B)_0$ that is full  in ${\rm M}_n(B)_0$ and for which $A \grcong e{\rm M}_n(B)e$.  An application of one direction of Theorem \ref{generalalgBGRthm} then finishes the proof.

Before pursuing the proof of Theorem~\ref{hgeTheorem} we examine the relationship between graded equivalence and homogeneously graded equivalence.  To begin, the following example shows that homogeneously graded equivalence is strictly stronger than graded equivalence.

\begin{example}\label{genothgeingeneral}
Let $E$ be the graph with two vertices $v$ and $w$, an edge from $v$ to $w$, and an edge from $w$ to $v$.  Note that $v L_k(E) v$ is a graded subring of $L_k(E)$ and $v$ is a full homogeneous idempotent of $L_k(E)$.  By \cite[Example~2.3.2]{RoozbehBook}, $v L_k(E)v$ is graded equivalent to $L_k(E)$.

We claim that $v L_k(E)v$ is not homogeneously graded  equivalent to $L_k(E)$.  Assume to the contrary that $v L_k(E)v$ is homogeneously graded  equivalent to $L_k(E)$.  By Theorem~\ref{hgeTheorem}, ${\rm M}_\infty(vL_k(E)v)$ and ${\rm M}_\infty( L_k(E))$ would then be graded isomorphic.  But this is a contradiction since straightforward computations yield that  ${\rm M}_\infty(vL_k(E)v)_1 = \{ 0 \}$ and ${\rm M}_\infty(L_k(E))_1 \neq \{ 0 \}$.  Hence, $v L_k(E)v$ is not homogeneously graded  equivalent to $L_k(E)$.

We note that while $v \in (L_K(E))_0$ is a full idempotent in $L_K(E)$, it is not a full idempotent in the subring $ (L_K(E))_0$ (as $w \notin (L_K(E))_0 v (L_k(E))_0$).   We also note, with Proposition \ref{geimplieshge} below in mind, that $L_k(E)$ is strongly graded by \cite[Theorem 1.6.13]{RoozbehBook}, but $vL_k(E)v$ is not (since $(vL_k(E)v)_1 = (vL_k(E)v)_{-1} = \{0\}$).
\end{example}

The previous example notwithstanding, there is an important situation in which a graded equivalence between $Gr-A$ and $Gr-B$ automatically implies that $A$ and $B$ are homogeneously graded equivalent.   Recall that the $\Gamma$-graded ring $A$ is called {\it strongly graded} when $A_\gamma A_\delta = A_{\gamma + \delta}$ for all $\gamma, \delta \in \Gamma$.   

\begin{proposition}\label{geimplieshge}
Suppose $A$ and $B$ are strongly  $\Gamma$-graded rings.   Then $A$ and $B$ are graded equivalent if and only if $A$ and $B$ are homogeneously graded equivalent.   
\end{proposition}

\begin{proof}
By Dade's Theorem (see, e.g., \cite[Theorem 1.5.1]{RoozbehBook}), the strongly graded hypothesis gives that  the functors
$$ - \otimes_{A_0} A : Mod-A_0 \to Gr-A  \ \ \  (\mbox{resp., }   - \otimes_{B_0} B : Mod-B_0 \to Gr-B)$$ are equivalences, with inverse functors given by the restriction functors  
$$( - )_0 : Gr-A \to Mod-A_0 \ \ \  (\mbox{resp., } ( - )_0 : Gr-B \to Mod-B_0).$$
 So if $\psi: Gr-A \to Gr-B$ is a graded equivalence, then the functor $\eta:  Mod-A_0 \to Mod-B_0$ given by  
 $$\eta := \  ( - )_0  \circ   \psi \circ ( - \otimes_{A_0}A) $$
is an equivalence and has the desired commutativity property with $\psi$. 
\end{proof}

\begin{corollary}\label{strongLpa}
Let $E$ and $F$ be finite graphs that each have no sinks.  Then $L_K(E)$ is graded equivalent to $L_K(F)$ if and only if $L_K(E)$ is homogeneously graded equivalent to $L_K(F)$.    
\end{corollary}

\begin{proof}
By \cite[Theorem 1.6.13]{RoozbehBook}, for any finite graph $E$, the Leavitt path algebra $L_K(E)$ is strongly graded if (and only if) $E$ has no sinks.   The result then follows from Proposition~\ref{geimplieshge}.  
\end{proof}

\begin{remark}\label{directsumsaredifferent}
For any ring $R$ and $n\in \N$ let $R^{(n)}$ denote the direct sum of $n$ copies of $R$ (viewed as a ring).  It is well known that for any field $K$, if $m, n\in \N$ for which $K^{(m)}$ is Morita equivalent to $K^{(n)}$, then $m=n$.  
\end{remark}

\begin{example}\label{specificgraphsnotgradedequiv}
For the graphs

$ $
\[
\xymatrix{
E: &  \bullet  \ar@(lu,ru)[] & \qquad \qquad \text{and} \qquad \qquad & F: & \bullet \ar@/^/[r] &  \bullet \ar@/^/[l].
}
\]
the Leavitt path algebra $L_k(E)$ is not homogeneously graded equivalent to the Leavitt path algebra $L_k(F)$ because
\[
L_k(E)_0 \cong k \qquad \qquad \text{and} \qquad \qquad L_k(F)_0 \cong k \oplus k,
\]
so that by Remark \ref{directsumsaredifferent} we get that $L_k(E)_0$ is not Morita equivalent to $L_k(F)_0$.   Thus by Corollary~\ref{strongLpa}, $L_k(E)$ is not graded equivalent to $L_k(F)$ either.
\end{example}

Along the way to establishing Theorem \ref{hgeTheorem}, we invoke a specific result established in the first three paragraphs of the proof  of \cite[Theorem 2.3.7]{RoozbehBook}.    

\begin{proposition}\label{Roozbehprojectivesprop}
Suppose there exists a graded category equivalence $ \tau: Gr-A \to Gr-B$.  If $P := \tau (A)$, then 
$A \grcong \operatorname{End}_B(P).$ 
\end{proposition}

    \begin{remark}
     We note that the proof of  \cite[Theorem 2.3.7]{RoozbehBook}  (and thus the validity of Proposition ~\ref{equivalentversionhge} and of Proposition~\ref{Roozbehprojectivesprop})  is highly nontrivial.  For instance, within the proof of \cite[Theorem 2.3.7]{RoozbehBook} it is necessary to use the nontrivial fact that if a graded right $B$-module is finitely generated projective  (respectively, a generator)  as an object in $Gr-B$, then it is also necessarily a finitely generated projective (respectively, a generator) as an object in $Mod-B$.   This is not at all obvious.
    \end{remark}

\noindent \emph{Proof of $(HG4) \implies (HG1)$ in Theorem~\ref{hgeTheorem}.} 
We begin by showing that the following diagram commutes whenever $e$ is a full idempotent in $A_0$. 
\[   \xymatrix{Mod-A_0 \ar[d]_{- \otimes_{A_0}A}\ar[rr]^-{- \otimes_{A_0}A_0e} & &Mod-eA_0e \ar[d]^{- \otimes_{eA_0e}eAe}\\
Gr-A\ar[rr]^-{- \otimes_AAe} & & Gr-eAe  }   \]
Rephrased, we show that 
$$ (M\otimes_{A_0}A_0e) \otimes_{eA_0e}eAe \grcong (M\otimes_{A_0}A) \otimes_A Ae$$
for every right $A_0$-module $M$.   
To do so, 
we define the  maps $\Gamma_1, \Gamma_2, \Delta_1$, and $\Delta_2$ as follows:   
$$\Gamma_1:  (M\otimes_{A_0}A_0e) \otimes_{eA_0e}eAe \to  M\otimes_{A_0}(A_0e \otimes_{eA_0e}eAe) \ ,  \ \ \  \mbox{via}    \ \ (m \otimes x) \otimes y \mapsto m \otimes (x\otimes y),$$
$$\Gamma_2:  (M\otimes_{A_0}A) \otimes_{A}Ae \to   M\otimes_{A_0} (A \otimes_{A}Ae) \ , \ \ \  \mbox{via} \ \  (m\otimes a) \otimes be \mapsto m\otimes (a \otimes be),$$
$$\Delta_1: M\otimes_{A_0}(A_0e \otimes_{eA_0e}eAe) \to M\otimes_{A_0}Ae \ , \ \ \ \mbox{via} \ \  m \otimes (x \otimes y) \mapsto m\otimes xy \ \ \mbox{ for } x\in A_0e, y\in eAe, \ \ \ \ \  \text{ and } 
$$
$$\Delta_2: M\otimes_{A_0} (A \otimes_{A}Ae) \to  M\otimes_{A_0}Ae  \ , \ \ \ \mbox{via} \ \  m\otimes (a \otimes be) \mapsto m\otimes abe \ \mbox{ for } a , b \in A.
$$

\noindent
Each of these maps is  an isomorphism of abelian groups that preserves the appropriate module structures (on both sides).  This is easy to see for $\Gamma_1^{-1}$ and $\Gamma_2^{-1}$, as they are just the inverses of the associativity maps.  To establish that both $\Delta_1$ and $\Delta_2$ are isomorphisms, we explicitly construct their inverses.       Using the fullness hypothesis, we have that  $1_A = 1_{A_0} = \sum_{i=1}^n a_ieb_i $ for some $a_i, b_i \in A_0$.   Define
$$\Theta_1 : M\otimes_{A_0}Ae \to M\otimes_{A_0}(A_0e \otimes_{eA_0e}eAe) \ , \ \  \ \mbox{via} \ \  m\otimes \alpha e \mapsto m \otimes (\sum_{i=1}^n a_ie\otimes eb_i \alpha e) \ \mbox{ for } \alpha \in A$$
and
$$ \Theta_2:   M\otimes_{A_0}Ae \to M\otimes_{A_0} (A \otimes_{A}Ae) \ , \ \ \  \mbox{via} \ \  m\otimes \alpha e \mapsto  m \otimes (\sum_{i=1}^n a_ie\otimes eb_i \alpha e) \ \mbox{ for } \alpha \in A.$$
At first glance, both $\Theta_1$ and $\Theta_2$ seem to produce the same outputs, but for $\Theta_1$ the output is viewed in 
$M\otimes_{A_0}(A_0e \otimes_{eA_0e}eAe) $, while for $\Theta_2$ the output is viewed in $M\otimes_{A_0} (A \otimes_{A}Ae)$. 

We shall establish that $\Theta_1 \circ \Delta_1$ is the identity on $M\otimes_{A_0}(A_0e\otimes_{eA_0e}eAe)$ by showing that $\Theta_1 \circ \Delta_1$ is the identity on pure tensors.   Let  $m\in M, x\in A_0$, and $\alpha \in A$.   
  Then 
\begin{eqnarray*}
 \ \ \ \ \ \ \ \  \Theta_1 \circ \Delta_1 \Big( m\otimes_{A_0} (xe \otimes_{eA_0e} e\alpha e) \Big) & = & \Theta_1(m\otimes_{A_0} xe\alpha e)  \\
& = & m \otimes_{A_0}  \sum_{i=1}^n \Big(a_ie  \otimes_{eA_0e} eb_ixe\alpha e \Big) \\
& = &   \sum_{i=1}^n \Big( m \otimes_{A_0} (a_ie \otimes_{eA_0e} eb_ixe\alpha e) \Big)  \\
& = &   \sum_{i=1}^n \Big( m \otimes_{A_0} (a_i eb_ix e \otimes_{eA_0e} e\alpha e )\Big)  \ \ \   (\mbox{since}  \  b_i , x \in A_0  \Rightarrow  eb_ixe \in eA_0e ) \\
& = &  m \otimes_{A_0} \Big( (\sum_{i=1}^n a_ieb_i ) xe  \otimes_{eA_0e} e\alpha e \Big)\\
& = & m\otimes_{A_0} (xe \otimes_{eA_0e} e\alpha e ) \ \hspace{.8in}  (\mbox{since} \sum_{i=1}^n a_ieb_i = 1_{A_0})   \ .  \\
  \end{eqnarray*}
\noindent
A similar computation  shows that  $\Delta_1 \circ \Theta_1$ is the identity on $M\otimes_{A_0}Ae$, which yields that $\Theta_1 = \Delta_1^{-1}$.  In a similar manner, one can establish $\Theta_2 = \Delta_2^{-1}.$
Now define the abelian group isomorphisms
$$\Phi:  (M\otimes_{A_0}A_0e) \otimes_{eA_0e}eAe \  \to \   (M\otimes_{A_0}A) \otimes_{A}Ae \ , \ \ \ \Phi := \Gamma_2^{-1} \circ \Delta_2^{-1} \circ \Delta_1 \circ \Gamma_1,$$
and
$$\Psi:  (M\otimes_{A_0}A) \otimes_{A}Ae  \ \to \  (M\otimes_{A_0}A_0e) \otimes_{eA_0e}eAe \ , \ \ \ \Psi := \Gamma_1^{-1} \circ \Delta_1^{-1} \circ \Delta_2 \circ \Gamma_2.$$

\noindent
The fact that both $\Phi$ and $\Psi$ are graded module homomorphisms follows directly from the grading imposed on the right $eAe$-modules   $(M\otimes_{A_0}A_0e) \otimes_{eA_0e}eAe$ and $(M\otimes_{A_0}A) \otimes_{A}Ae$, respectively.

Thus we have established that the diagram above commutes for all objects in $Mod - A_0$.  For ease of application in the next part of the proof, we explicitly state the following fact, which is a consequence of the commutativity of the diagram above for any unital rings $R,S$.  
 
\smallskip
 
\noindent {\bf Fact}:   Suppose  $e \in R_0$ is full in $R_0$, and let $S = eRe$.    Then for every right $R_0$-module $M$, we have $ (M\otimes_{R_0}R_0e) \otimes_{S_0}S  \cong   (M\otimes_{R_0}R) \otimes_{R}Re$ as graded right $S$-modules.   

\smallskip

\noindent We shall use the Fact,  together with the proof of the (HG4) $\implies$ (HG2) implication of Theorem~\ref{generalalgBGRthm} to establish the (HG4) $\implies$ (HG1) implication of Theorem~\ref{hgeTheorem}.

Assuming (HG4), let $\varphi: {\rm M}_\infty(A) \to {\rm M}_\infty(B)$ denote the indicated graded isomorphism.    
 Define $e = \varphi(e^A_{1,1}) \in {\rm M}_\infty(B)$.   We proceed in steps:
 
\noindent (1)   Because the grading on all matrix rings is standard, we have,  via the restriction $\varphi |$ of $\varphi$ to $e_{1,1}{\rm M}_\infty(A)e_{1,1}$,  that $A \grcong e{\rm M}_\infty(B)e \grcong e{\rm M}_n(B)e$ where $n \in \mathbb{N}$ is such that $e \in {\rm M}_n(B)$.       So via  $\varphi |$ we get that   the following diagram commutes (on objects):
\[   \xymatrix{Mod-A_0 \ar[d]_{- \otimes_{A_0}A}\ar[r]^-{  \eta} &Mod-e{\rm M}_n(B)_0e \ar[d]^{- \otimes_{e{\rm M}_n(B)_0e}e{\rm M}_n(B)e}\\
Gr-A\ar[r]^-{    \psi} &Gr-e{\rm M}_n(B)e  }   \]
\noindent
(Note that we have not used the Fact in this step.)

\smallskip

 \noindent (2)   As established  in the proof of the (HG4) $\implies$ (HG2) implication of Theorem~\ref{generalalgBGRthm}, we have that  $e$ is full in ${\rm M}_n(B)_0$ for an appropriate $n \in \mathbb{N}$.   We then apply the Fact with $R := {\rm M}_n(B)$ and $S := e{\rm M}_n(B)e$ to deduce the following diagram commutes on objects:
 
 \[   \xymatrix{Mod-e{\rm M}_n(B)_0e        \ar[d]_{- \otimes_{e{\rm M}_n(B)_0e}e{\rm M}_n(B)e}       \ar[r]^-{\xi} &Mod-{\rm M}_n(B)_0            \ar[d]^{- \otimes_{{\rm M}_n(B)_0}{\rm M}_n(B)}          \\
Gr-e{\rm M}_n(B)e\ar[r]^-{\tau} &Gr-{\rm M}_n(B)  }   \]
 
\noindent (Note that we have written the diagram in the opposite order from the one that is used in the Fact; it is really the inverses of these functions that come into play here.  We denote those inverses by $\xi$ and $\tau$, respectively.)

\smallskip

\noindent (3)   We see that $e_{1,1} = e^B_{1,1}$ is full in ${\rm M}_n(B)_0$. Using $R:={\rm M}_n(B)$ and $S := e^B_{1,1}{\rm M}_n(B)e^B_{1,1} = B$, the Fact implies that the following diagram commutes:

\[   \xymatrix{Mod-{\rm M}_n(B)_0 \ar[d]_{- \otimes_{{\rm M}_n(B)_0}{\rm M}_n(B)}    \ar[rrr]^-{- \otimes_{{\rm M}_n(B)_0}{\rm M}_n(B)_0e_{1,1}}   &  &    &  Mod-e_{1,1}{\rm M}_n(B)_0e_{1,1}    \ \cong Mod-B_0 \ar[d]^{- \otimes_{e_{1,1}{\rm M}_n(B)_0e_{1,1}}e_{1,1}{\rm M}_n(B)e_{1,1}}  \\
Gr-{\rm M}_n(B)\ar[rrr]^-{- \otimes_{{\rm M}_n(B)}{\rm M}_n(B)e_{1,1}} & &    &    Gr-e_{1,1}{\rm M}_n(B)e_{1,1} \ \cong Gr-B   }   \]

\bigskip

Upon composing the three diagrams above, we deduce that the diagram

\[   \xymatrix{Mod-A_0 \ar[d]_{- \otimes_{A_0}A}\ar[r]^-{\eta} &Mod-B_0 \ar[d]^{- \otimes_{B_0}B}\\
Gr-A\ar[r]^-{\psi} &Gr-B }   \]

\noindent
commutes on objects, where  $\eta$ and $\psi$ are each built as the appropriate compositions of  functors.      Since each of the constituent functors in these compositions is an equivalence (respectively, graded equivalence), we conclude that $\eta$ (respectively, $\psi$) is an equivalence (respectively, graded equivalence) of categories, thereby completing the proof of the (HG4) $\implies$ (HG1) implication of Theorem~\ref{hgeTheorem}.  
\qed

\bigskip

\noindent \emph{Proof of $(HG1) \implies (HG4)$ in Theorem~\ref{hgeTheorem}.}   We assume there exists a category equivalence $\varphi_1 : Mod-A_0 \to Mod-B_0$ and a graded category equivalence $\varphi_2: Gr-A \to Gr-B$ for which the following diagram commutes on objects:

\[   \xymatrix{Mod-A_0 \ar[d]_{- \otimes_{A_0}A}\ar[r]^-{\varphi_1} &Mod-B_0 \ar[d]^{- \otimes_{B_0}B}\\
Gr-A\ar[r]^-{\phi_2} &Gr-B }   \]

\noindent By Proposition \ref{Roozbehprojectivesprop} we deduce $A \grcong \operatorname{End}_B(P)$,  where $P := \varphi_2(A)$.     (Note:   To conclude this, we only use the existence of a graded equivalence between $Gr-A$ and $Gr-B$;  i.e., we only need the information coming from the bottom row of the diagram).

   We show  $A_0 \otimes_{A_0} A \grcong A$ as graded right $A$-modules, via the following argument.       The isomorphism $A_0 \otimes_{A_0} A \to A$ given by $x \otimes a \mapsto xa$, is standard.  The fact that this isomorphism preserves the grading is clear, because the grading on $A_0 \otimes_{A_0} A$ is defined by setting $(A_0 \otimes_{A_0} A)_\gamma := A_0 \otimes_{A_0} (A_\gamma)$ for each $\gamma \in \Gamma$.

  So together  with the hypothesis that the above diagram commutes we obtain
      $$ P := \varphi_2(A) \grcong \varphi_2(A_0\otimes_{A_0}A) \grcong \varphi_1(A_0) \otimes_{B_0} B \ \quad\mbox{as graded right $B$-modules}.$$

      Let $x_1, x_2, \dots , x_n$ be a set of generators of $\varphi_1(A_0)$ as a right $B_0$-module.    Such a finite set exists because $A_0$ is finitely generated as a right $A_0$-module (since $A_0$ is unital) and $\varphi_1$ is a category equivalence.    We conclude that $x_1 \otimes 1_B$, $x_2\otimes 1_B$, $\dots$, $x_n\otimes 1_B$ is a set of generators for $\varphi_1(A_0) \otimes_{B_0} B$ as a right $B$-module.  Moreover, by the definition of the grading on the tensor product,   each of these generators is in the $0$-component of $\varphi_1(A_0) \otimes_{B_0} B$ . Since $P \grcong \varphi_1(A_0) \otimes_{B_0} B$,  we conclude that there is a set of generators $p_1, p_2, \dots , p_n$ of $P$ as a right $B$-module, for which each element of the generating set is in the $0$-component of $P$.

      Because (as noted above) $P$ is a (finitely generated) projective module in the category $Mod-B$, the map $\tau:  B^n \to P$ taking $(b_1, b_2, \dots , b_n) \mapsto \sum_{i=1}^n p_ib_i$ necessarily splits; let $\sigma$ denote an inverse map, so that $ \tau \circ \sigma $ is the identity on $P$.   Then there is a right $B$-module $Q$ with $P \oplus Q \cong B^n$ as right $B$-modules.   Because each $p_i$ is in $P_0$,  the map $\tau: (B(0))^n \to P$ is easily seen to be a graded homomorphism (where $B(0)$ is the graded module $B_B$ with no translation applied).  We also see that $\sigma$ is graded; that is,
      $$P \oplus Q \grcong (B(0))^n.$$
      
 \noindent
  By the argument given in \cite[page 37]{RoozbehBook}, we deduce that this graded isomorphism of right $B$-modules yields a graded isomorphism of graded rings, and
$$\operatorname{End}_B ( (B(0) ) ^n ) \grcong {\rm M}_n(B),$$ 
where ${\rm M}_n(B)$ has the standard grading.

      Let  $f$ denote the endomorphism $f \in \operatorname{End}_B ( (B(0) ) ^n )$ that has $f |_P = 1_P$ and $f |_Q = 0$, and let $e$  be the  idempotent in ${\rm M}_n(B)_0 = {\rm M}_n(B_0)$ that is the image of $f$ under the displayed isomorphism. Then   
     \begin{eqnarray*}
     A  & \grcong &\operatorname{End}_B(P) \qquad  \mbox{(from above)} \\
      & \grcong & f \operatorname{End}_B ( (B(0) ) ^n ) f   \\ 
      & \grcong &  e {\rm M}_n(B) e. 
  \end{eqnarray*}
Since $P_0 \cong \varphi_1(A_0)$ as right $B_0$-modules,  we conclude $e$ is full in ${\rm M}_n(B)_0$.

 Finally, we apply the (HG2) $\implies$ (HG4) implication of the  Algebraic Stabilization Theorem (i.e., Theorem \ref{generalalgBGRthm}), to conclude ${\rm M}_\infty(A) \grcong {\rm M}_\infty(B)$ as desired.
 \qed

 \bigskip

We end this section by identifying various situations in which pairs of Leavitt path algebras are homogeneously graded equivalent.    

Let $E$ be a graph.  For a non-negative integer $n$, we let $E^n$ denote the paths of length $n$ (where paths of length zero are the vertices of the graph).  The set of all paths of $E$ will be denoted by ${\rm Path}(E)$; thus ${\rm Path}(E) = \bigcup_{n=0}^\infty E^n$.

For each positive integer $n$, we let $G_n$ denote the graph with $n+1$ vertices $z_0, z_1, \ldots, z_n$,  and exactly one edge from $z_i$ to $z_{i-1}$.
$$G_n :  \qquad  \xymatrix{ \bullet^{z_n} \ar[r] &  \bullet^{z_{n-1}}  \ar[r] & \cdots  \ar[r]& \bullet^{z_1} \ar[r] & \bullet^{z_0}} \ .$$

\begin{corollary}\label{cor:linegraph}
Let $E$ and $F$ be finite acyclic graphs.  Suppose $E$ has exactly one sink $v$ and $F$ has exactly one sink $w$.  Then $L_k(E)$ is homogeneously graded  equivalent to $L_k(F)$ if and only if 
$$\max\{ {\rm length}(\mu) : \mu \in {\rm Path}(E), r(\mu) = v \} \ = \ \max \{ {\rm length} ( \nu ) : \nu \in {\rm Path}(F), r(\nu) =  w \}.$$
\end{corollary}    

\begin{proof}
Suppose $ \max\{ \operatorname{length}(\mu) : \mu \in {\rm Path}(E), r(\mu)=v \} \  =  \ \max \{ \operatorname{length}( \nu ) : \nu  \in {\rm Path}(F) \text{ and } r(\nu)=w \}$; denote this common integer by $n$.  If $n = 0$, then both $E$ and $F$ are graphs with one vertex and no edges,  so  we trivially have that $L_k(E) \cong k$ is graded isomorphic (and also homogeneously graded equivalent) to $L_k(F) \cong k$ .  Assume that $n \geq 1$.  Let $E_1$ be the graph obtained from $G_n$ by attaching a source to $z_j$ whenever 
\[
n^E_j = | \{ \mu : \mu \in E^{j+1} \text{ and } r(\mu)=v \} | \geq 2\]
and adding $n^E_j -1$ edges from this source to $z_{j-1}$.  Let $F_1$ be the graph obtained from $G_n$ by attaching a source to $z_{j}$ whenever
\[
n^F_j = | \{ \mu : \mu \in F^{j+1} \text{ and } r(\mu)=w \} | \geq 2
\]
and adding $n^E_j -1$ edges from this source to $z_{j}$.  By \cite[Theorem 4.13]{Hazrat}, $L_k(E)$ is graded isomorphic (and hence homogeneously graded equivalent) to $L_k(E_1)$, and $L_k(F)$ is graded isomorphic (and hence homogeneously graded equivalent) to $L_k(F_1)$.  

We next show that $L_k(E_1)$ is homogeneously graded equivalent to $L_k(G_n)$.  We see that $G_n$ is a subgraph of $E_1$ in such a way that $L_k(G_n)$ is graded isomorphic to $p L_k(E_1)p$   where $p = \sum_{ i = 0}^n z_i$.  We claim that $p$ is full in $L_k(E_1)_0$.  Let $v$ be a source of $E_1$ that is not $z_n$.  Let $\mu$ be a path in $E_1$ starting at $v$ and ending at $v_0$, and let $m = \operatorname{length}(\mu)$.  By construction of $E_1$, there exists a path $\nu$ from $z_{ m}$ to $z_0$, and this path has length $m$.  Thus, $x=\mu \nu^* \in L_k(E_1)_0$ such that $xx^* = \mu \mu^*$ and $x^*x = \nu \nu^* = z_{ m }$.  Thus $\mu \mu^* \in (L_k(E_1))_0 p (L_k( E_1))_0$.  Since $E_1$ is a finite acyclic graph with exactly one sink,
\[
v = \sum_{ \substack{ \mu \in {\rm Path}(E) \\ s(\mu)=v , r(\mu)=z_0} } \mu \mu^* \in (L_k(E_1))_0 p (L_k( E_1))_0.
\]
Consequently, $ (L_k(E_1))_0 p (L_k( E_1))_0=  (L_k(E_1))_0$.   By Theorem~\ref{generalalgBGRthm} and Theorem~\ref{hgeTheorem}, we deduce $L_k(E_1)$ and $pL_k(E_1)p \grcong L_k(G_n)$ are homogeneously graded equivalent.  
A similar argument shows that $L_k(G_n)$ is homogeneously graded equivalent to $L_k(F_1)$.  Hence $L_k(E)$ and $L_k(F)$ are homogeneously graded equivalent.

Conversely, assume $L_k(E)$ is homogeneously graded equivalent to $L_k(F)$.  Arguing as above, $L_k(E)$ is homogeneously graded  equivalent to $L_k(G_m)$, and $L_k(F)$ is homogeneously graded  equivalent to $L_k(G_n)$, where
\[m = \max\{ \operatorname{length}(\mu) : \mu \in {\rm Path}(E), r(\mu)=v \}  \quad \text{and} \quad
n = \max\{ \operatorname{length}(\mu) : \mu \in {\rm Path}(F), r(\mu) = w \}. 
\]
Since $L_k(E)$ is homogeneously graded  equivalent to $L_k(F)$, we conclude $L_k(G_m)$ is homogeneously graded  equivalent to $L_k(G_n)$.   Hence there is an equivalence of categories between $Mod-L_k(G_m)_0$ and $Mod-L_k(G_n)_0$.  Note that for all $\mu, \nu \in {\rm Path}(G_m)$, if  $\mu \nu^*$ is nonzero in $L_k(G_m)$ and $\operatorname{length}(\mu) = \operatorname{length}( \nu )$, then $\mu =\nu$.   Therefore,
\begin{align*}
L_k(G_m)_0 &= \mathrm{span} \{ \mu \nu^* :  \mu, \nu \in {\rm Path}(G_m) ,  \operatorname{length}(\mu) = \operatorname{length}(\nu) \} = \mathrm{span} \{ \mu \mu^* : \mu \in {\rm Path}(G_m) \} \\
&= \mathrm{span} \{ s(\mu) : \mu \in {\rm Path}(G_m) \} = \mathrm{span} \{  z_i : 0 \leq i \leq m \} \cong \bigoplus_{i=0}^m k.
\end{align*}  
Similarly, $L_k(G_n)_0 \cong  \bigoplus_{i=0}^n k$.   Hence there is an equivalence of categories between $Mod-\bigoplus_{i=0}^m k$ and $Mod-\bigoplus_{i=0}^n k$, and since $k$ is a field, this implies that $m=n$ (see Remark \ref{directsumsaredifferent}). 
\end{proof}

\begin{example}\label{specificgraphshge} 
Our results allows us to deduce when the Leavitt path algebras of various acyclic graphs are graded equivalent as well as homogeneously graded equivalent.

\begin{enumerate}
\item Consider the graphs
\[
\xymatrix{
E : & \bullet & \qquad \text{and}  \qquad &
F : & \bullet \ar[r] & \bullet
}
\]
Then $L_k(E)$ is graded equivalent to $L_k(F)$ by \cite[Example 2.3.2]{RoozbehBook} and \cite[Theorem 4.13]{Hazrat} but it follows from Corollary~\ref{cor:linegraph} that $L_k(E)$ is not homogeneously graded equivalent $L_k(F)$.

\item Consider the graphs
\[
\xymatrix{
G :  & \bullet  \ar[r] & \bullet \ar[r] & \bullet & \qquad \text{and} \qquad & H : & \bullet \ar[r] & \bullet \ar[r] & \bullet \\
  &  & && & &	& \bullet \ar[u] \ar[r] & \bullet \ar[u]   \\
 &   & && & & &	& \bullet \ar[u] 
}
\]
It follows from Corollary~\ref{cor:linegraph} that $L_k(G)$ is homogeneously graded equivalent $L_k(H)$.
\end{enumerate}

\end{example}

\section{The Homogeneously Graded Version of the Extended Morita Theorem}

In this section we establish our main result, which is stated in Theorem~\ref{MainTheorem}.  This theorem establishes the equivalence of the appropriate  analogues, in the context of homogeneously graded equivalence,  of the statements appearing in the Extended Morita Theorem presented in the Introduction.    In the previous two sections we have established some of the main components of this result.    To establish the remaining connections, we introduce additional key ideas in the following definitions.

\begin{definition}\label{sequencegradingdef}    Let $m\in \N \cup \{\infty\}$.   Let $(\delta_n)_{1 \leq n \leq m}$ (for $m$ finite) or $(\delta_n)_{n \in \NN}$ (for $m = \infty$) be a sequence in $\Gamma$.   
Let $A$ be a $\Gamma$-graded ring.  We let  ${\rm M}_m(A)[ (\delta_n) ]$  denote the ring ${\rm M}_m(A)$ with the following $\Gamma$-grading:  for each $\lambda \in \Gamma$, 
\[
({\rm M}_m(A)[(\delta_n)])_\lambda := ( A_{ \lambda+\delta_j - \delta_i} )_{i,j} .
\]
In particular, suppose that $(\delta_n)$ is a constant sequence,  e.g., $\delta_n = 0$ for all $n\in \NN$.  Then ${\rm M}_m(A)[(\delta_n)]$ is precisely ${\rm M}_m(A)$ with the standard grading.  
\end{definition}

\begin{remark}\label{cornergrading}
Let $m\in \N \cup \{\infty\}$.   Let $(\delta_n)_{1 \leq n \leq m}$ (for $m$ finite) or $(\delta_n)_{n \in \NN}$ (for $m = \infty$) be a sequence in $\Gamma$. Then, by construction, for any $\Gamma$-graded ring $A$, the corner ring $e_{1,1}({\rm M}_m(A)[(\delta_n)])e_{1,1}$ is graded isomorphic to $A$.  
\end{remark}

\begin{definition}
Let $A$ and $B$ be unital $\Gamma$-graded rings.  A \emph{graded (surjective) Morita context}  $(A, B, M, N, \psi, \phi)$ between $A$ and $B$ consists of graded bimodules $_AM_B$ and $_BN_A$, a surjective graded $A$-homomorphism $\psi \colon M\otimes_B N \to A$, and a surjective graded $B$-homomorphism $\phi \colon N \otimes_A M \to B$ satisfying 
\[
\phi( n \otimes m)n' = n \psi(m \otimes n') \quad \text{and} \quad m' \phi( n \otimes m) = \psi(m' \otimes n) m
\]
for all $m, m' \in M$ and $n , n' \in N$.
\end{definition}

\begin{definition}
Let $A$ and $B$ be unital $\Gamma$-graded rings.  A \emph{homogeneously graded (surjective) Morita context}  $(A, B, M, N, \psi, \phi)$ between $A$ and $B$ consists of a graded (surjective) Morita context $(A, B, M, N, \psi, \phi)$ between $A$ and $B$ such that $(A_0, B_0, M_0, N_0 , \psi_{ M_0 \otimes N_0} , \phi_{ N_0 \otimes M_0})$ is a (surjective) Morita context, where $\psi_{ M_0 \otimes N_0}$ is the composition of the maps $M_0 \otimes N_0 \to (M\otimes N)_0$ and $\psi\vert_{ (M \otimes N)_0}$, and $\phi_{ N_0 \otimes M_0}$ is the composition of the maps $N_0 \otimes N_0 \to (N\otimes M)_0$ and $\psi\vert_{ (N \otimes M)_0}$. 
\end{definition}

\begin{lemma}\label{lem:moritacontext}
Let $A$, $B$, and $C$ be unital $\Gamma$-graded rings.  Suppose there exists a homogeneously graded (surjective) Morita context $(A, B, M, N, \psi, \phi)$ between $A$ and $B$ and there exists a homogeneously graded (surjective) Morita context $(B, C, X, Y, \beta, \alpha)$ between $B$ and $C$.  Then $(A, C, M \otimes X , Y \otimes N, \nu, \mu )$ is a homogeneously graded (surjective) Morita context between $A$ and $C$, where $\mu ( (m \otimes x) \otimes (y \otimes n) ) = \psi( m \otimes (\beta(x \otimes y ) n))$ and $\nu ( (y \otimes n) \otimes (m \otimes x) ) = \alpha(y \otimes (\phi( n \otimes m )x))$.
\end{lemma}

\begin{proof}
A computation shows that $\mu$ is an $A$-homomorphism.  We shall show $\mu$ is a graded $A$-homomorphism.  Let $\gamma \in \Gamma$.  Assume that $m \in M_{ \gamma_1}$, $x \in X_{ \gamma_2}$, $y \in Y_{\gamma_3}$, and $n \in N_{\gamma_4}$ such that $\gamma_1 + \gamma_2 +\gamma_3 + \gamma_4 = \gamma$.  Since $\beta$ is a graded $B$-homomorphism, $\beta( x \otimes y) \in B_{ \gamma_2+\gamma_3}$.  Hence, $(\beta(x \otimes y ) n) \in N_{ \gamma_2 + \gamma_3 + \gamma_4 }$ which implies that $m \otimes (\beta(x \otimes y ) n)) \in (M \otimes N)_{ \gamma_1 + \gamma_2 + \gamma_3 + \gamma_4 } = (M \otimes N)_\gamma$.  Since $\psi$ is a graded $A$-homomorphism, we have $\mu ( (m \otimes x) \otimes (y \otimes n) ) = \psi( m \otimes (\beta(x \otimes y ) n)) \in A_\gamma$.  Thus $\mu$ is a graded $A$-homomorphism, because $((M\otimes X) \otimes (Y \otimes N))_\gamma$ is spanned by elements of the form $(m \otimes x) \otimes (y \otimes n)$, where $m \in M_{ \gamma_1}$, $x \in X_{ \gamma_2}$, $y \in Y_{\gamma_3}$, and $n \in N_{\gamma_4}$ satisfying $\gamma_1 + \gamma_2 +\gamma_3 + \gamma_4 = \gamma$.  Similarly, $\nu$ is a graded $C$-homomorphism.

Since $\beta \colon X \otimes Y \to B$ is a surjective $B$-homomorphism, and since $N$ is an $A -B$-bimodule, every element in $N$ is in the linear span of $\beta( x \otimes y) n$, where $x \in X$, $y \in Y$, and $n\in N$.  Using this fact and the fact that $\psi \colon M \otimes N \to A$ is a surjective $A$-homomorphism, we deduce that $A$ is in the linear span of $\mu( (m \otimes x) \otimes (y \otimes n))= \psi( m \otimes (\beta(x \otimes y ) n))$ where $m \in M$, $x \in X$, $y \in Y$, and $n \in N$.  Therefore, $\mu \colon (M\otimes X) \otimes (Y \otimes N) \to A$ is a surjective $A$-homomorphism.  A similar argument shows that $\nu \colon (Y \otimes N) \otimes (M\otimes X) \to C$ is a surjective $C$-homomorphism.

To show $\mu_{ (M \otimes X)_0 \otimes (Y \otimes N)_0 } \colon (M \otimes X)_0 \otimes (Y \otimes N)_0 \to A_0$ is a surjective homomorphism, note that $B_0 N_0 = N_0$ and  $(M\otimes X)_0 \otimes (Y \otimes N)_0$ contains elements of the form $(m \otimes x) \otimes (y \otimes n)$ where $m \in M_0$, $x \in X_0$, $y \in Y_0$, and $n \in N_0$.  Also observe that, by the definition of a homogeneously graded (surjective) Morita context, $\beta_{X_0 \otimes Y_0 } \colon X_0 \otimes Y_0 \to B_0$ and $\psi_{ M_0 \otimes N_0 } \colon M_0 \otimes N_0 \to A_0$ are surjective homomorphisms.  We use an argument similar to that used in the above paragraph to deduce that $\mu_{ (M \otimes X)_0 \otimes (Y \otimes N)_0 } \colon (M \otimes X)_0 \otimes (Y \otimes N)_0 \to A_0$ is a surjective homomorphism.  Likewise, a similar argument gives $\nu_{ (Y \otimes N)_0 \otimes (M\otimes X)_0 } : (Y \otimes N)_0 \otimes (M\otimes X)_0 \to C_0$ is a surjective homomorphism.  Hence $(A, C, M \otimes X , Y \otimes N, \nu, \mu )$ is a homogeneously graded (surjective) Morita context between $A$ and $C$.
\end{proof}

\begin{lemma}\label{lem:fullcorner-moritacontext}
Let $A$ be a unital $\Gamma$-graded ring, and let $e$ be a full idempotent in $A_0$.  Then there are graded surjective homomorphisms $\psi \colon eA \otimes Ae \to e Ae$ and $\phi \colon Ae \otimes eA \to A$ such that $\psi ( x \otimes y) = xy$, $\phi(y \otimes x)=yx$, and $(e Ae, A, eA , Ae, \psi , \phi )$ is a homogeneously graded (surjective) Morita context between $eAe$ and $A$. 
\end{lemma}

\begin{proof}
Since $(eA)_\gamma = eA_{\gamma}$ and $(Ae)_\gamma= A_\gamma e$, we conclude that $\psi$ is a graded $eAe$-homomorphism and $\phi$ is a graded $A$-homomorphism such that $\psi_{ (eA)_0 \otimes (Ae)_0 } \colon (eA)_0 \otimes (Ae)_0 \to eA_0e$ is surjective and the image of $\phi_{ (Ae)_0 \otimes (eA)_0 } \colon (Ae)_0 \otimes (eA)_0 \to A_0$ is $A_0 e A_0$.  Since $e$ is a full idempotent in $A_0$, we have $A_0 e A_0 = A_0$.  Thus, $\phi_{ (Ae)_0 \otimes (eA)_0 }$ is a surjective $A_0$-homomorphism.

Since $A_0$ is a unital subalgebra of $A$, $A_0$ is full in $A$.  Hence, $e$ is full in $A$.  A similar argument as above shows that $\psi$ is a surjective $eAe$-homomorphism and $\phi$ is a surjective $A$-homomorphism.  Thus $(e Ae, A, eA , Ae, \psi , \phi )$ is a homogeneously graded (surjective) Morita context between $eAe$ and $A$. 
\end{proof}

We now have both the precise terminology and the mathematical machinery to establish our main result:  The Homogeneously Graded Version of the Extended Morita Theorem, as described in the Introduction.

\begin{theorem}\label{MainTheorem}
Let $A$ and $B$ be unital $\Gamma$-graded  rings.  Give all matrix rings the standard grading.    Then the following are equivalent.

\begin{itemize}

\item[(HG1)\, ] $A$ and $B$ are homogeneously graded equivalent. 

\item[(HG2)\, ] There exist $n\in \N$ and a (homogeneous) idempotent $e \in {\rm M}_n(B)_0$ that is full in ${\rm M}_n(B)_0$   for which the rings $A$ and $ e{\rm M}_n(B)e$ are graded isomorphic. 
 
\item[(HG3)\, ] There exists a homogeneously graded (surjective) Morita context between $A$ and $B$. 
 
\item[(HG4)\, ] ${\rm M}_\infty(A) $ is graded isomorphic to ${\rm M}_\infty(B)$.  
 
\item[(HG4*)] For any sequence $(\delta_n)_{ n \in \NN }$ with the property that the set $\{ k : \delta_k = \delta_n \}$ is infinite for each $n \in \NN$, we have ${\rm M}_\infty(A)[(\delta_n)]$ is graded isomorphic to ${\rm M}_\infty(B)[(\delta_n)]$.
\end{itemize} 
 
 \end{theorem}

 \begin{proof}
The equivalence of (HG2) and (HG4) is precisely The Algebraic Stabilization Theorem (Theorem~\ref{generalalgBGRthm}).  The equivalence of (HG1) and (HG4) is precisely Theorem~\ref{hgeTheorem}.   To establish Theorem \ref{MainTheorem}, we shall show that (HG2) implies (HG3), that (HG3) implies (HG1), and that (HG4) and (HG4*) are equivalent.

We first show that (HG2) implies (HG3).  Assume there exist $n\in \N$ and a (homogeneous) idempotent $e \in {\rm M}_n(B)_0$ that is full in ${\rm M}_n(B)_0$ and for which the rings $A$ and $ e{\rm M}_n(B)e$ are graded isomorphic.  Let $\{ e_{i,j} \}$ be the system of matrix units in $M_n(B)_0 = M_n(B_0)$, where the $(i,j)$-entry of $e_{i,j}$ is $1_B$ and all other entries are zero.  Then $e_{1,1}$ is a full (homogeneous) idempotent in ${\rm M}_n(B)_0$ such that $B$ is graded isomorphic to $e_{1,1} {\rm M}_n(B) e_{1,1}$.  By Lemma~\ref{lem:moritacontext} and Lemma~\ref{lem:fullcorner-moritacontext}, there exists a homogeneously graded (surjective) Morita context between $e {\rm M}_n(B) e$ and $e_{1,1} {\rm M}_n(B) e_{1,1}$ which then induces a homogeneously graded (surjective) Morita context between $A$ and $B$ due to the facts that $A$ is graded isomorphic to $e {\rm M}_n(B) e$ and $B$ is graded isomorphic to $e_{1,1} {\rm M}_n(B) e_{1,1}$.  Thus (HG3) holds.

Next we show that (HG3) implies (HG1).  Assume that there exists a homogeneously graded (surjective) Morita context,  $(A, B, M, N, \psi, \phi)$, between $A$ and $B$.   Consider the Morita ring $L= \begin{pmatrix} A & M \\ N & B \end{pmatrix}$, which is $\Gamma$-graded by setting $L_\gamma = \begin{pmatrix} A_\gamma & M_\gamma \\ N_\gamma & B_\gamma \end{pmatrix}$ for each $\gamma \in \Gamma$, with addition defined entrywise and multiplication given by
\[
\begin{pmatrix}
a & m \\ n & b 
\end{pmatrix}
\begin{pmatrix}
a' & m' \\
n' & b' 
\end{pmatrix}
=
\begin{pmatrix}
aa' + \psi(m' \otimes n') & am' + m b' \\
n a' + b n' & \phi(n \otimes m') + bb'
\end{pmatrix}. 
\]
Since $L_0 = \begin{bmatrix} A_0 & M_0 \\ N_0 & B_0 \end{bmatrix}$ and $(A_0, B_0, M_0, N_0 , \psi_{ M_0 \otimes N_0} , \phi_{ N_0 \otimes M_0})$ is a (surjective) Morita context, the elements $p = \begin{pmatrix} 1_A & 0 \\ 0 & 0 \end{pmatrix}$ and $q= \begin{pmatrix} 0 & 0 \\ 0 & 1_B \end{pmatrix}$ are full idempotents in $L_0$.   Hence, by the equivalence of (HG1) and (HG2), $p L p$ is homogeneously graded equivalent to $L$, and $qLq$ is homogeneously graded equivalent to $L$.  By the previously-noted transitivity property of homogeneously graded equivalence,  we get that  $pLp$ is homogeneously graded equivalent to $qLq$.  Since $A$ is graded isomorphic to $pLp$, and since $B$ is graded isomorphic to $qLq$  (both statements use the definition of the grading on $L$), we conclude $A$ is homogeneously graded equivalent $B$.  Thus (HG1) holds.

The equivalence of (HG4) and (HG4*) will be established below in Proposition~\ref{HG4andHG4*}.  
 \end{proof}

\section{The Graded Version of the Extended Morita Theorem}

In this the final section of the article we extend a result of \cite{RoozbehBook} to establish The Graded Version of the Extended Morita Theorem described in the Introduction.    We then conclude the article by establishing the equivalence of statements (HG4) and (HG4*) of Theorem~\ref{MainTheorem}.  

\begin{lemma}\label{lem:11corner}
Let $A$ be a $\Gamma$-graded ring, let $\{ e_{i,j} \}$ be a system of matrix units for ${\rm M}_\infty(A)$, and let $(\delta_n)_{ n \in \NN }$ be  a sequence in $\Gamma$.  Then $A$ is graded isomorphic to $e_{1,1} {\rm M}_\infty(A) [ (\delta_n) ] e_{1,1}$.
\end{lemma}

\begin{proof}
Let $a\otimes e_{1,1}$ 
be the element in ${\rm M}_\infty(A)$ that has $a$ in the $(1,1)$-entry and 0 in all other entries.  Then $\phi \colon A \to e_{1,1} {\rm M}_\infty(A)[(\delta_n)] e_{1,1}$,  defined by setting  $\phi(a)=a\otimes e_{1,1}$ for each $a\in A$,  is a ring isomorphism.  Since
$$
(e_{1,1} {\rm M}_\infty(A) [ (\delta_n) ] e_{1,1})_\lambda = e_{1,1} {\rm M}_\infty(A) [ (\delta_n ) ] e_{1,1} \cap {\rm M}_\infty(A)[(\delta_n)]_\lambda = \begin{pmatrix} A_\lambda & \mathbf{0} \\ \mathbf{0} & \mathbf{0}  \end{pmatrix},
$$
$\phi$ is a graded isomorphism.
\end{proof}

\begin{theorem}\label{thm:gr-me-stablization}
Let $R$ and $S$ be $\Gamma$-graded rings.  Then $R$ is graded equivalent to $S$ if and only if there exists a sequence $(\gamma_n)_{n \in \NN }$ in $\Gamma$ such that ${\rm M}_\infty(R)[(0)]$ is graded isomorphic to ${\rm M}_\infty(S) [ (\gamma_n) ]$.
\end{theorem}

\begin{proof}
Suppose there exists a sequence $(\gamma_n)_{n \in \NN }$ in $\Gamma$ such that ${\rm M}_\infty(R)[(0)]$ is graded isomorphic to ${\rm M}_\infty(S) [ (\delta_n) ]$.  Then 
\begin{align*}
R &\cong_\mathrm{gr} e_{1,1} {\rm M}_\infty(R)[(0)] e_{1,1} \qquad \text{(by Lemma~\ref{lem:11corner})} \\
&\cong_\mathrm{gr} f {\rm M}_\infty(S)[(\gamma_n)] f \quad \text{for some full homogeneous idempotent } f \text{ in } {\rm M}_\infty(S)[(\gamma_n)] \\
&\cong_\mathrm{gr} \overline{f} {\rm M}_n(S)[(\gamma_n)] \overline{f} \quad \text{for some full homogeneous idempotent } \overline{f} \text{ in } {\rm M}_n(S) [ \gamma_1, \cdots , \gamma_n ].
\end{align*}
(Here $\overline{f}$ is simply $f \in {\rm M}_\infty(S)$, viewed in a finite matrix ring of appropriate size.) Since $\overline{f}$ is a full homogeneous idempotent in ${\rm M}_n(S) [ \gamma_1, \cdots , \gamma_n ]$, by \cite[Theorem 2.3.8]{Hazrat} we conclude that $\overline{f} {\rm M}_n(S)[ \gamma_1, \cdots , \gamma_n ] \overline{f}$ is graded Morita equivalent to ${\rm M}_n(S) [ \gamma_1, \cdots , \gamma_n ]$.  Again using \cite[Theorem 2.3.8]{Hazrat},  we have $e_{1,1} {\rm M}_n(S)[ \gamma_1, \cdots , \gamma_n ] e_{1,1}$ is graded Morita equivalent to ${\rm M}_n(S) [ \gamma_1, \cdots , \gamma_n ]$.  Since $S$ is graded isomorphic to $e_{1,1}{\rm M}_n(S) [ \gamma_1, \cdots , \gamma_n ] e_{1,1}$ (see Remark \ref{cornergrading}),  we conclude $R$ is graded equivalent to $S$.

Conversely, suppose $R$ is graded equivalent to $S$.  By \cite[Theorem 2.3.8]{Hazrat} there exist  $n \in \NN$, elements $\delta_1, \cdots, \delta_n \in \Gamma$, and a full homogeneous idempotent $e$ of the graded ring ${\rm M}_n(S)[\delta_1, \cdots, \delta_n ]$ such that $R$ is graded isomorphic to the corner $e {\rm M}_n(S)[\delta_1, \cdots , \delta_n ] e$.  Since $1_R$ is full in $R_0$, the idempotent $e$ is full in ${\rm M}_n(S)[\delta_1, \cdots, \delta_n ]_0$.  Thus, by The Algebraic Stabilization Theorem (Theorem~\ref{generalalgBGRthm}), we have
\[
{\rm M}_\infty(e {\rm M}_n(S)[ \delta_1, \cdots, \delta_n] e) [(0)] \cong_\mathrm{gr} {\rm M}_\infty({\rm M}_n(S)[ \delta_1, \cdots, \delta_n] ) [(0)].
\]

Define $(\gamma_n)_{n \in \NN}$ by $\gamma_k = \gamma_r$ where $1 \leq r \leq n$ and $k \equiv r \ (\operatorname{mod} n)$.  Define $\phi \colon {\rm M}_\infty({\rm M}_n(S)[\delta_1, \cdots, \delta_n ] ) [(0)] \to {\rm M}_\infty(S) [ ( \gamma_k) ]$ by 
\begin{align*}
\phi( (A_{i,j} )_{i,j} ) = \begin{pmatrix} A_{1,1} & A_{1,2} & \cdots  \\
A_{2,1} & A_{2,2} & \cdots  \\
\vdots & \vdots & \ddots \end{pmatrix}.
\end{align*}
A computation shows that $\phi$ is a graded isomorphism.  Hence, using the previously displayed isomorphism, we get 
$$ {\rm M}_\infty(R)[(0)] \cong_\mathrm{gr} {\rm M}_\infty(e {\rm M}_n(S)[ \delta_1, \cdots, \delta_n] e) [(0)] 
\cong_\mathrm{gr} {\rm M}_\infty({\rm M}_n(S)[ \delta_1, \cdots, \delta_n] ) [(0)]
\cong_\mathrm{gr} {\rm M}_\infty(S) [ ( \gamma_k) ],
$$
which establishes the result.
\end{proof}

\begin{theorem}\label{GradedVersionExtMoritaThm}   Let $A$ and $B$ be $\Gamma$-graded unital rings.  Then the following are equivalent.  
\begin{itemize}
\item[(GM1)] $A$ and $B$ are graded equivalent.     

\item[(GM2)] There exist $n\in \N$ and an  idempotent $e \in {\rm M}_n(B)$ that is full in  ${\rm M}_n(B)$ and a sequence $(\gamma_m)_{1\leq m \leq n} \in \Gamma$ for which the rings $A$ and $ e{\rm M}_n(B)[(\gamma_m)]e$ are graded isomorphic.  
  
\item[(GM3)] There exists a graded (surjective) Morita context between $A$ and $B$. 
 
\item[(GM4)] There exists a sequence $(\gamma_m)_{m \in \NN }$ in $\Gamma$ such that ${\rm M}_\infty(A)[(0)]$ is graded isomorphic to ${\rm M}_\infty(B) [ (\gamma_m) ]$.
\end{itemize}

\end{theorem}

\begin{proof}

The equivalence of (GM1) and (GM2) is established in \cite[Theorem 2.3.8]{RoozbehBook}.  The equivalence (GM1) and (GM3) is argued mutatis mutandis as in the argument that (HG1) is equivalent to (HG3) in the proof of Theorem~\ref{MainTheorem} by simply disregarding the adjective ``homogeneously".  The equivalence of (GM1)  and (GM4) is Theorem~\ref{thm:gr-me-stablization}.  
\end{proof}

We note that the equivalence of (GM1) and (GM3) can also be established by arguing exactly as in the proof that (MT1) is equivalent to (MT3),  by simply inserting the adjective ``graded".  

We conclude with an analysis of various gradings on the infinite matrix ring ${\rm M}_\infty(A)$ for a $\Gamma$-graded ring $A$.   This leads us to a proof of the equivalence of statements (HG4) and (HG4*) in Theorem~\ref{MainTheorem}.

\begin{lemma}\label{lem:permutation}
Let $R$ be a unital $\Gamma$-graded ring, let $(\delta_n)_{n \in \NN}$ be a sequence in $\Gamma$, and let $\pi \colon \NN \to \NN$ be a bijection.  Then ${\rm M}_\infty(R) [ (\delta_n) ]$ is graded isomorphic to ${\rm M}_\infty(R)[ ( \delta_{\pi(n)} ) ]$.
\end{lemma}

\begin{proof}
Define $\phi \colon {\rm M}_\infty(R) [ (\delta_n) ] \to {\rm M}_\infty(R)[ ( \delta_{\pi(n)} ) ]$ by $
\phi( (a_{i, j}) ) = ( a_{ \pi(i), \pi(j) } )$.
Let $A = (a_{i,j } )$ and $B = (b_{i,j} )$ be matrices in ${\rm M}_\infty(R) [ (\delta_n) ]$.  Then 
\begin{align*}
\phi(AB) &= \phi \left( \left( \sum_{ l=1}^\infty a_{i,l} b_{l,j}  \right) \right)
        = \left( \sum_{ l =1}^\infty a_{ \pi(i) , l } b_{ l , \pi(j) } \right) 
        = \left( \sum_{ l =1}^\infty a_{ \pi(i) , \pi(l) } b_{ \pi(l) , \pi(j) } \right)  \\
        &= ( a_{ \pi(i), \pi(j) } ) ( b_{ \pi(i), \pi(j) } )
        = \phi(A)\phi(B) 
\end{align*}        
and
$$\phi(A+B) = \phi \left( \left( a_{i,j} + b_{i,j} \right)\right)
        = ( a_{\pi(i), \pi(j) } + b_{ \pi(i), \pi(j) } )
        = \phi(A)+\phi(B),$$
so that $\phi$ is a ring homormophism.  Suppose $(a_{i,j} ) \in {\rm M}_\infty(R) [ (\delta_n) ]_\lambda$.  Then $a_{i,j} \in A_{\lambda + \delta_j - \delta_i}$.  Therefore, $a_{ \pi(i) , \pi(j) } \in A_{\lambda + \delta_{\pi(j)} - \delta_{\pi(i)} }$, which implies that $( a_{ \pi(i), \pi(j) } ) \in {\rm M}_\infty(R)[ ( \delta_{\pi(n)} ) ]_\lambda$.  Hence $\phi$ is a graded isomorphism.
\end{proof}

\begin{lemma}\label{lem:double-stablization}
Let $A$ be a $\Gamma$-graded ring, and let $( \delta_n )_{n \in \NN}$ be a sequence in $\Gamma$.  Define $(\gamma_n)$ as follows.  Write $\NN = \bigsqcup_{ k = 1}^\infty S_k$ with $|S_k| = \infty$ for all $k \in \NN$.  Define $\gamma_n := \delta_k$ for each $n \in S_k$.  Then 
$${\rm M}_\infty (A) [(\gamma_n )] \mbox{ is graded isomorphic to }{\rm M}_\infty ( {\rm M}_\infty(A) [(0)] ) [ (\delta_n) ].$$

Consequently, if for each $n \in \NN$, the set $\{ k : \delta_k = \delta_n \}$ is infinite, then ${\rm M}_\infty ( A) [ (\delta_n) ]$ is graded isomorphic to ${\rm M}_\infty ( {\rm M}_\infty(A) [(0)] ) [ (\delta_n) ]$.
\end{lemma}

\begin{proof}
For each $k \in \NN$, chose a bijection $\pi_k \colon \NN \to S_k$.  (Such a bijection exists since $S_k$ is an infinite subset of $\NN$.)  For a ring $R$ and for each $i, j \in \NN$ and each $x \in R$, let $
x \otimes e_{i,j}$ denote the element in ${\rm M}_\infty(R)$ that has $x$ in the $(i,j)$-entry and zero in all other entries.  

Note that every element of ${\rm M}_\infty ( {\rm M}_\infty(A) [(0)] ) [ (\delta_n) ]$ can be uniquely written as a finite linear combination of elements of the form 
$(a \otimes e_{i,j}) \otimes e_{k,l}$
 for $a \in A$ and $i,j, k, l \in \NN$.  Therefore, the mapping
\[
(a \otimes e_{i,j}) \otimes e_{k,l} \mapsto a \otimes e_{ \pi_{k}(i), \pi_l(j) }
\]
extends to a linear function $\phi \colon {\rm M}_\infty ( {\rm M}_\infty(A) ) \to {\rm M}_\infty (A)$.  Note that
$$
((a \otimes e_{i,j}) \otimes e_{k,l}) \cdot ((b \otimes e_{i',j'}) \otimes e_{k',l'})
\  = \  \delta_{l,k'} ((a\otimes e_{i,j})(b\otimes e_{i',j'}) \otimes e_{k,l'} )
\ =  \
 \delta_{l, k'} \delta_{j,i'} ((ab \otimes e_{i,j'}) \otimes e_{k,l'})$$
 and
$$
(a \otimes e_{ \pi_{k}(i), \pi_l(j) }) (b \otimes e_{ \pi_{k'}(i'), \pi_{l'}(j') }) \ = \ \delta_{ \pi_l(j) , \pi_{k'}(i') } (ab \otimes e_{ \pi_k(i) , \pi_{l'}(j') }) \ = \  \delta_{l, k'} \delta_{j,i'} (ab \otimes e_{ \pi_k(i) , \pi_{l'}(j') }),
$$
where the last equality follows from the facts that each $\pi_k$ is a bijection and the ranges of $\pi_k$ and $\pi_{l'}$ are disjoint when $k \neq l'$.  Consequently, $\phi$ is a ring homomorphism.

Let $a \in A$, and let $r,s \in \NN$.  Then $r \in S_k$ and $s \in S_l$ for some $k, l \in \NN$, beacuse $\NN = \bigsqcup_{ k = 1}^\infty S_k$.  Therefore, since $\pi_k, \pi_l$ are bijections, there exist $i, j \in \NN$ such that $\pi_k(i)=r$ and $\pi_l(j)=s$.  This implies 
\[
\phi( (a \otimes e_{i,j}) \otimes e_{l,k}) \ = \ a\otimes e_{r,s}.
\]
Consequently, $\phi$ is surjective.  To show injectivity, suppose $x \in {\rm M}_\infty({\rm M}_\infty(A))$ with $\phi(x)=0$.  Writing  $x = \sum_{ i,j,k,l \in \NN} (a_{i,j,k,l} \otimes e_{i,j}) \otimes e_{k,l}$
 with only finitely many of the $a_{i,j,k,l}$ terms nonzero, 
\[
0 = \phi(x) = \sum_{ i,j,k,l \in \NN} a_{i,j,k,l} \otimes e_{ \pi_{k}(i), \pi_l(j) }.
\]
Since the ranges of the $\pi_k$ maps are disjoint and each $\pi_k$ is a bijection, 
\[
\sum_{ i,j,k,l \in \NN} a_{i,j,k,l} \otimes e_{ \pi_{k}(i), \pi_l(j) } 
= 0
\]
for all $i,j,k,l \in \NN$.  Consequently, $\phi$ is injective, and $\phi$ is a ring isomorphism.

Let $\lambda \in \Gamma$, and suppose $x \in {\rm M}_\infty ( {\rm M}_\infty (A)[(0)] )[ (\delta_n) ]_\lambda$.  As above, write $x =\sum_{ i,j,k,l \in \NN} (a_{i,j,k,l} \otimes e_{i,j}) \otimes e_{k,l}$ with only finitely many of the $a_{i,j,k,l}$ being nonzero.  For each $i,j, k,l \in \NN$, we have
\[
a_{i,j, k,l} \otimes e_{i,j}
 \in {\rm M}_\infty(A)[(0)]_{\lambda + \delta_l - \delta_k} = {\rm M}_\infty(A_{\lambda + \delta_l - \delta_k}).
\]
Thus for each $i,j, k,l \in \NN$, we have $a_{i,j, k, l} \in A_{\lambda + \delta_l - \delta_k}$. Since for each $i, k \in \NN$, we have $\gamma_{\pi_k(i) } = \delta_k$, we conclude
\[
 A_{\lambda + \delta_l - \delta_k} = A_{\lambda + \gamma_{ \pi_l(j) } - \gamma_{ \pi_k(i) } }
\]
for all $i,j, k,l \in \NN$.  Hence for each $i, j, k, l \in \NN$, we have $
a_{i,j,k, l} \otimes e_{\pi_k(i), \pi_l(j)} 
 \in A_{\lambda + \gamma_{ \pi_l(j) } - \gamma_{ \pi_k(i) } }$.  Consequently,
\[
\phi ( x ) = \sum_{ i , j , k,l \in \NN} 
a_{i,j,k, l}  \otimes e_{\pi_k(i), \pi_l(j)}
 \in {\rm M}_\infty(A)[(\gamma_n)]_\lambda.
\]
Thus $\phi$ is a graded isomorphism.

For the final part of the lemma, suppose that for each $n \in \NN$, the set $\{ k : \delta_k = \delta_n \}$ is infinite. Let $\pi \colon \NN \to \NN$ be a bijection with $( \delta_{\pi(n)} )_{ n \in \NN} = (\gamma_n)_{ n \in \NN }$.  By Lemma~\ref{lem:permutation}, ${\rm M}_\infty(A)[ (\gamma_n) ]$ is graded isomorphic to ${\rm M}_\infty (A)[(\delta_n)]$.  Since ${\rm M}_\infty( {\rm M}_\infty(A)[(0)] ) [ (\delta_n) ]$ is graded isomorphic to ${\rm M}_\infty(A)[ (\gamma_n) ]$, by composing the graded isomorphisms we obtain our desired result that 
${\rm M}_\infty ( A) [ (\delta_n) ] \mbox{ is graded isomorphic to }{\rm M}_\infty ( {\rm M}_\infty(A) [(0)] ) [ (\delta_n) ].$
\end{proof}

\begin{proposition}\label{HG4andHG4*}
Let $A$ and $B$ be $\Gamma$-graded unital rings.   Then the following are equivalent.
\begin{itemize}
\item[(HG4)\, ]   ${\rm M}_\infty(A) $ is graded isomorphic to ${\rm M}_\infty(B)$ in the standard grading.  
\item[(HG4*)]   For any sequence $(\delta_n)_{ n \in \NN }$ with the property that the set $\{ k : \delta_k = \delta_n \}$ is infinite for each $n \in \NN$, we have ${\rm M}_\infty(A)[(\delta_n)]$ is graded isomorphic to ${\rm M}_\infty(B)[(\delta_n)]$.
\end{itemize}

\end{proposition}

\begin{proof}
We remind the reader that the standard grading on ${\rm M}_\infty(A)$ (respectively, ${\rm M}_\infty(B)$) is precisely the grading represented by the notation  ${\rm M}_\infty(A)[(0)] $ (respectively, ${\rm M}_\infty(B)[(0)] )$, where $(0)$ is the constant sequence with each entry equal to $0$.

 Suppose (HG4) holds.  Then ${\rm M}_\infty(A) = {\rm M}_\infty(A)[(0)]$ is graded isomorphic to ${\rm M}_\infty(B) = {\rm M}_\infty(B)[(0)]$.  Let $(\delta_n)_{ n \in \NN }$ be a sequence with $\{ k : \delta_k = \delta_n \}$ infinite for each $n \in \NN$.  Then
 \begin{align*}
{\rm M}_\infty(A)[(\delta_n)] &\cong_\mathrm{gr} {\rm M}_\infty( {\rm M}_\infty(A)[(0)] ) [(\delta_n)] \qquad \text{ (by Lemma~\ref{lem:double-stablization})} \\
&\cong_\mathrm{gr} {\rm M}_\infty( {\rm M}_\infty(B)[(0)] ) [(\delta_n)] \\
&\cong_\mathrm{gr} {\rm M}_\infty(B)[(\delta_n)] \qquad \qquad \text{ (by Lemma~\ref{lem:double-stablization}).}
\end{align*}
Hence (HG4) implies (HG4*).

To establish (HG4*) implies (HG4*), simply note that the constant sequence $(\delta_n)_{n \in \NN} = (0)$ trivially has the property that $\{ k : \delta_k = \delta_n\} = \NN$ is infinite for each $n \in \NN$.  Thus (HG4*) implies ${\rm M}_\infty(A)[(0)]$ is graded isomorphic to ${\rm M}_\infty(B)[(0)]$, and (HG4) holds.
\end{proof}

\end{document}